\documentclass[12pt]{amsart}
\usepackage{amssymb,latexsym}
\usepackage{enumerate}
\usepackage{amscd}

\theoremstyle{plain}
\newtheorem{theorem}{Theorem}
\newtheorem{corollary}[theorem]{Corollary}
\newtheorem{lemma}[theorem]{Lemma}
\newtheorem{proposition}[theorem]{Proposition}
\theoremstyle{definition}

\newtheorem{example}[theorem]{Example}

\theoremstyle{remark}
\newtheorem*{remark}{Remark}

\begin{document}

% One author
\title[Regular Coverings]{Regular Coverings in
Filter and Ideal Lattices}
\author{William H. Rowan}
\address{POB 20161 \\
         Oakland, California 94620-0161}
\email{whrowan@member.ams.org}
%\thanks{thanks}
% End one author

% Some input statements - I am not sure where they should
% go.
%\input mymathdefs3

\def\pair#1#2{\langle #1, #2\rangle}
\def\triple#1#2#3{\langle #1, #2, #3\rangle}
\def\triplem#1#2#3{M(#1,#2,#3)}
\def\triplemp#1{M(#1)}
\def\semitimes{\ltimes}
\def\rsemitimes{\rtimes}
\def\Mof#1{\Cal M(A)}
\def\pp#1#2{{\roman P[#1,#2]}}
\def\inclusion#1#2{{\scriptstyle[#1\hookrightarrow #2]}}

\newcommand{\N}{{\Bbb N}}
\newcommand{\Z}{{\Bbb Z}}
\newcommand{\Q}{{\Bbb Q}}
\newcommand{\R}{{\Bbb R}}
\newcommand{\A}{{\Bbb A}}
%\def\F{{\roman F}}
%\def\U{{\roman U}}

%category names
\def\Set{{\text{\bf Set}}}
\def\Mon{{\text{\bf Mon}}}
\def\Grp{{\text{\bf Grp}}}
\def\Rng{{\text{\bf Rng}}}
\def\CRng{{\text{\bf CRng}}}
\def\Ab{{\text{\bf Ab}}}
\def\Pnt{{\text{\bf Pnt}}}
\def\Alg{{\text{\bf Alg}}}
\def\Smgrp{{\text{\bf Smgrp}}}
\def\Ov{{\text{\bf Ov}}}
\def\Pre{{\text{\bf Pre}}}
\def\Lat{{\text{\bf Lat}}}
\def\LatX{{\text{\bf Lat-X}}}
\def\Modlat{{\text{\bf Modlat}}}
\def\Typ{{\text{\bf Typ}}}
\def\catClo{{\text{\bf Clo}}}
\def\Mod{{\text{\bf Mod}}}
\def\Bimod{{\text{\bf Bimod}}}
\def\Lie{{\text{\bf Lie}}}
\def\Nonass{{\text{\bf Nonass}}}
\def\Rep{{\text{\bf Rep}}}
\def\Rngd{{\text{\bf Rngd}}}
\def\Algd{{\text{\bf Algd}}}
\def\Cat{{\text{\bf Cat}}}
\def\Bimod{{\text{\bf Bimod}}}
\def\Abhp{{\text{\bf AbHp}}}
\def\Clnd{{\text{\bf Clnd}}}
\def\Sktch{{\text{\bf Sktch}}}
\def\Top{{\text{\bf Top}}}
\def\Vect{{\text{\bf Vect}}}
\def\BPSet{{\text{\bf BPSet}}}
\def\Abmod{{\text{\bf Abmod}}}
\def\Unif{{\text{\bf Unif}}}
\def\HUnif{{\text{\bf HUnif}}}
\def\CHUnif{{\text{\bf CHUnif}}}
\def\CU{{\text{\bf CU}}}
\def\HCU{{\text{\bf HCU}}}
\def\CHCU{{\text{\bf CHCU}}}

%other category-theoretic stuff
\def\natur{{\overset{\scriptscriptstyle\bullet}\to\to}}
%(was \text{\bf .} instead)

% Operator names:
\newcommand{\Cg}{\operatorname{Cg}}
\newcommand{\Con}{\operatorname{Con}}
\newcommand{\Cov}{\operatorname{Cov}}
\newcommand{\Div}{\operatorname{Div}}
\newcommand{\End}{\operatorname{End}}
\newcommand{\Eqv}{\operatorname{Eqv}}
\newcommand{\Ext}{\operatorname{Ext}}
\newcommand{\Fg}{\operatorname{Fg}}
\newcommand{\Hom}{\operatorname{Hom}}
\newcommand{\Ig}{\operatorname{Ig}}
\newcommand{\Mg}{\operatorname{Mg}}
\newcommand{\Pol}{\operatorname{Pol}}
\newcommand{\OpSemiUnif}{\operatorname{SemiUnif}}
\newcommand{\Sub}{\operatorname{Sub}}
\newcommand{\Sg}{\operatorname{Sg}}
\newcommand{\Term}{\operatorname{Term}}
\newcommand{\coker}{\operatorname{coker}}
\newcommand{\Ke}{\operatorname{Ke}}
\newcommand{\Co}{\operatorname{Co}}
\newcommand{\kerc}{\operatorname{kerc}}
\newcommand{\kerp}{\operatorname{kerp}}
\newcommand{\MyIm}{\operatorname{Im}}
\newcommand{\nat}{\operatorname{nat}}
\newcommand{\trcl}{\operatorname{tr.\ cl.}}
\newcommand{\typ}{\operatorname{typ}}
\newcommand{\ass}{\operatorname{ass}}
\newcommand{\card}{\operatorname{card}}
\newcommand{\Rab}{\operatorname{Rab}}
\newcommand{\Clo}{\operatorname{Clo}}
\newcommand{\id}{\operatorname{1}}
\newcommand{\Id}{\operatorname{1}}
\newcommand{\Idl}{\operatorname{Idl}}
\newcommand{\Fil}{\operatorname{Fil}}
\newcommand{\Der}{\operatorname{Der}}
\newcommand{\Sym}{\operatorname{Sym}}
\newcommand{\Res}{\operatorname{Res}}
\newcommand{\Ind}{\operatorname{Ind}}
\newcommand{\Pro}{\operatorname{Pro}}
\newcommand{\Ann}{\operatorname{Ann}}
\newcommand{\Ob}{\operatorname{Ob}}
\newcommand{\Ar}{\operatorname{Ar}}
\newcommand{\cod}{\operatorname{cod}}
\newcommand{\dom}{\operatorname{dom}}
\newcommand{\core}{\operatorname{core}}
\newcommand{\subcore}{\operatorname{subcore}}
\newcommand{\Fit}{\operatorname{Fit}}
\newcommand{\Frat}{\operatorname{Frat}}
\newcommand{\Lim}{\operatorname{Lim}}
\newcommand{\Colim}{\operatorname{Colim}}
\newcommand{\rR}{\operatorname{R}}
\newcommand{\rU}{\operatorname{U}}
\newcommand{\rClndD}{\operatorname{Clnd}^{\Bbb D}}
\newcommand{\rClnd}{\operatorname{Clnd}}
\newcommand{\rSktchD}{\operatorname{Sktch}^{\Bbb D}}
\newcommand{\rSktch}{\operatorname{Sktch}}
\newcommand{\Reg}{\operatorname{Reg}}
\newcommand{\DDown}{\operatorname{Ddn}}
\newcommand{\DUp}{\operatorname{Dup}}
\newcommand{\Gr}{\operatorname{Gr}}
\newcommand{\Step}{\operatorname{Step}}
\newcommand{\nullset}{\{\}}
\newcommand{\Ass}{\operatorname{Ass}}
\newcommand{\cov}{\operatorname{cov}}
\newcommand{\Aut}{\operatorname{Aut}}
\newcommand{\Nor}{\operatorname{Nor}}
\newcommand{\I}{\operatorname{I}}
\newcommand{\OpUnif}{\operatorname{Unif}}
\newcommand{\Succ}{\operatorname{Succ}}

\renewcommand{\subjclassname}{\textup{2000} Mathematics Subject
     Classification}

\def\Var{{\mathbf V}}
\def\Str{{\hbox{\bf str}}}
\def\Spl{{\hbox{\bf ret}}}
\def\Ret{{\hbox{\bf ret}}}

%special notations
\def\lsil{\lbrack\!\lbrack}
\def\rsil{\rbrack\!\rbrack}

%hyphenation stuff
\def\congruence{on\-gru\-ence\discretionary{-}{}{-}}
\def\conM/{c\congruence mod\-u\-lar}
\def\ConM/{C\congruence mod\-u\-lar}
\def\conD/{c\congruence dis\-trib\-u\-tive}
\def\ConD/{C\congruence dis\-trib\-u\-tive}
\def\conP/{c\congruence per\-mut\-a\-ble}
\def\ConP/{C\congruence per\-mut\-a\-ble}
\def\conMity/{\conM/\-i\-ty}
\def\ConMity/{\ConM/\-i\-ty}
\def\conDity/{c\congruence dis\-trib\-u\-tiv\-i\-ty}
\def\ConDity/{C\congruence dis\-trib\-u\-tiv\-i\-ty}
\def\conPity/{c\congruence per\-mut\-a\-bil\-i\-ty}
\def\ConPity/{C\congruence per\-mut\-a\-bil\-i\-ty}
\def\usprv/{un\-der\-ly\-ing-set-pre\-ser\-ving}

% some common latin abbreviations
\def\ie/{{i.e.}}
\def\Ie/{{I.e.}}
\def\eg/{{e.g.}}
\def\Eg/{{E.g.}}
\def\etc/{{etc.}}

% power set notation
\def\pow{{\Cal P}}

% arrows underneath formulas
\newdimen\mysubdimen
\newbox\mysubbox
\def\subright#1{\mathpalette\subrightx{#1}}
\def\subleft#1{\mathpalette\subleftx{#1}}
\def\subrightx#1#2{\subwhat#1{#2}\rightarrowfill}
\def\subleftx#1#2{\subwhat#1{#2}\leftarrowfill}

\def\subwhat#1#2#3{{
\setbox\mysubbox=\hbox{#3}
\mysubdimen=\wd\mysubbox
\setbox\mysubbox=\hbox{$#1#2$}
\ifnum\mysubdimen>\wd\mysubbox
\vtop{
\hbox to\mysubdimen{\hfil\box\mysubbox\hfil}
\nointerlineskip
\hbox{#3}}
\else
\mysubdimen=\wd\mysubbox
\vtop{
\box\mysubbox
\nointerlineskip
\hbox to\mysubdimen{#3}}
\fi
}}

\date{\today}

\begin{abstract}
The Dedekind--Birkhoff theorem for finite-height
modular lattices has previously been generalized to
complete modular lattices, using the theory of regular
coverings. In this paper, we investigate regular
coverings in lattices of filters and lattices of ideals,
and the regularization strategy--embedding the lattice
into its lattice of filters or lattice of ideals, thereby
possibly converting a covering which is not regular into a
covering which is regular. One application of the theory
is a generalization of the notion of chief factors, and
of the Jordan-Holder Theorem, to cases where the
modular lattice in question is of infinite height.
Another application is a formalization of the
notion of the steps in the proof of a theorem.
\end{abstract}
\maketitle

\section*{Introduction}

 The purpose of this paper is to further
develop and apply the theory of regular coverings in a
complete modular lattice, introduced in
\cite{R97}. The point of that theory is to generalize, to
complete modular lattices, some of the nice results
available for finite-height modular lattices.

For example, given a finite-height module $M$ over a ring
$R$ (i.e., a module having a finite-height lattice of
submodules), a \emph{composition series} of $M$ is a
(necessarily finite) sequence of submodules
\[\{\,0\,\}=M_0\subset M_1\subset\ldots\subset M_n=M\]
such that each quotient $M_i/M_{i-1}$ is simple. The
Jordan-Holder Theorem states that any two composition
series $\{\,M_i\,\}_{i=1}^n$, $\{\,M'_i\,\}_{i=1}^{n'}$
are the same length~$n=n'$ and the
quotients~$M_i/M_{i-1}$ can be paired with the
quotients~$M'_j/M'_{j-1}$ in such a way that
corresponding quotients are isomorphic.

The Jordan-Holder Theorem is an algebraic
version of the lattice-theoretic Dedekind-Birkhoff
Theorem. The lattice-theoretic correlates of the
composition series and the isomorphism of corresponding
quotients are maximal chains (maximal linearly-ordered
subsets) and projective equivalence of coverings.
The Dedekind-Birkhoff Theorem states that in any two
maximal chains in a finite-height modular lattice, the
lengths of the chains are the same and coverings in the
chains can be paired in such a way that corresponding
coverings are projectively equivalent.

 Unfortunately, the
Dedekind-Birkhoff Theorem can fail for infinite-height
modular lattices. For example, consider the modular
lattice
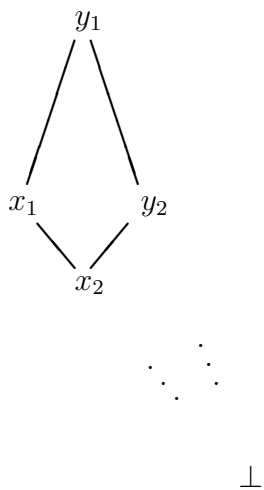
\[
\thicklines
\begin{picture}(135,200)
\put(50,190){$y_1$}
\put(25,120){$x_1$}
\put(75,120){$y_2$}
\put(50,90){$x_2$}
\put(77,60){$.$}
\put(82,54){$.$}
\put(87,48){$.$}
\put(96,68){$.$}
\put(99,61){$.$}
\put(102,54){$.$}
\put(112,15){$\bot$}
\put(56,184){\line(1,-3){18}}
\put(50,184){\line(-1,-3){18}}
\put(56,98){\line(5,6){14}}
\put(50,98){\line(-5,6){14}}
\end{picture}
\]
where $\bigwedge_ix_i=\bigwedge_iy_i=\bot$. There are two
maximal chains
\[C_1:\qquad\bot<\ldots<x_n<\ldots<x_2<x_1<y_1\]
and
\[C_2:\qquad\bot<\ldots<y_n<\ldots<y_2<y_1\]
such that the covering $x_1\prec y_1$ appears in $C_1$,
but there is no projectively equivalent covering in
$C_2$.

The theory of regular coverings was created to try to
remedy this situation. In the
language of that theory, the covering $x_1\prec y_1$ is
not
\emph{regular}, as defined in Section~\ref{S:RegCov}. If
it were regular, such behavior would be impossible
because of Theorem~\ref{T:MultInv}.

Now, note that if we embed the lattice $L$ into its
lattice of filters $\Fil L$, we obtain the lattice
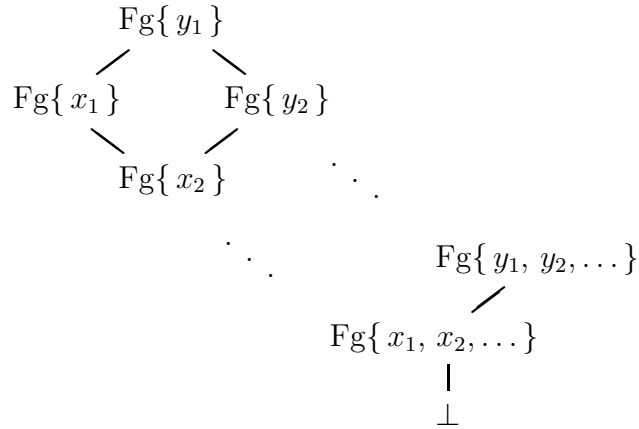
\[
\thicklines
\begin{picture}(220,180)
\put(50,160){$\Fg\{\,y_1\,\}$}
\put(10,130){$\Fg\{\,x_1\,\}$}
\put(90,130){$\Fg\{\,y_2\,\}$}
\put(50,100){$\Fg\{\,x_2\,\}$}
\put(130,40){$\Fg\{\,x_1,\,x_2,\ldots\,\}$}
\put(170,70){$\Fg\{\,y_1,\,y_2,\ldots\,\}$}
\put(170,10){$\bot$}
\put(130,108){$.$}
\put(138,102){$.$}
\put(146,96){$.$}
\put(90,78){$.$}
\put(98,72){$.$}
\put(106,66){$.$}
\put(42,143){\line(4,3){12}}
\put(98,143){\line(-4,3){12}}
\put(52,113){\line(-4,3){12}}
\put(83,113){\line(4,3){12}}
\put(184,53){\line(4,3){12}}
\put(175,23){\line(0,1){10}}
\end{picture}
\]
which contains new elements, $\Fg\{\,x_1,x_2,\ldots\,\}$
and $\Fg\{\,y_1,y_2,\ldots\,\}$, forming a covering
equivalent to $\Fg\{\,x\,\}\prec\Fg\{\,y\,\}$. Now any
two maximal chains in
$\Fil L$ have one covering equivalent to
$\Fg\{\,x\,\}\prec
\Fg\{\,y\,\}$. We have
\emph{regularized} the covering $x_1\prec y_1$ by
embedding $L$ into $\Fil L$. In this paper, we explore
this process and strategy of regularization further. Of
course, this involves studying coverings in $\Fil L$ (and
in its dual, the lattice $\Idl L$ of \emph{ideals} of
$L$) and trying to determine whether or not they are
regular.

We also examine questions of multiplicity, since in a
modular lattice that is not distributive, a maximal chain
can contain more than one covering projectively
equivalent to a given one.

 We give two applications. One application
is a generalization of the theory of chief factors of
an algebra having a modular congruence lattice. The
information supplied by these results is entirely
lattice-theoretic; we leave for another time the
algebraic correlates such as play roles in the
Jordan-Holder Theorem. The other application is a way of
defining the steps in the proof of a theorem. 
Any proof of the theorem from the same premises must
\emph{cover}, as we say, the same steps. Also,
from any set of instances of rules of inference which
covers the steps, a finite subset can be
selected and used to construct a proof.

After this introduction and a section of preliminaries,
this paper begins in Section~\ref{S:Mult} with some
definitions relating to multiplicities. Given a maximal
chain $C$ in the modular lattice $L$, and a covering
$x\prec y$, there is a corresponding multiplicity of
$x\prec y$ in $C$ which may vary with $C$, except in the
important case when $x\prec y$ is \emph{weakly regular}.
We also define notations for upper and lower bounds on
the multiplicity.

Section~\ref{S:RegCov} discusses the theory of
\emph{regular coverings}, which, due to a generalization
of the Dedekind-Birkhoff Theorem as given in \cite{R97},
are also weakly regular.

Section~\ref{S:Coverings} is a preliminary examination
of coverings in filter and ideal lattices. As our
strategy is to use regular coverings in such lattices
for various purposes, we must understand their basic
properties before attempting to determine whether or not
they are regular. In this section, among other things,
we classify filter and ideal coverings into three
categories: \emph{atomic}, \emph{quasi-atomic}, and
\emph{anomalous} coverings.

Section~\ref{S:Stability} gives proofs of the stability
of regularity, multiplicity when regular, and in some
cases the multiplicity upper bound, under the embedding
from
$L$ into $\Fil L$ or $\Idl L$.

Section~\ref{S:Relationship} proves a relationship
between the multiplicity upper bound of a covering
$x\prec y$ in $L$ and the multiplicity lower bound of
the corresponding filter or ideal covering. The
important consequence of this is that under appropriate
conditions, if the
multiplicity bound is infinite, then  any maximal chain
in the filter or ideal lattice will have an infinite
number of coverings equivalent to
$\Fg\{\,x\,\}\prec\Fg\{\,y\,\}$ or
$\Ig\{\,x\,\}\prec\Ig\{\,y\,\}$. This complements other
theorems which describe the behavior when the
multiplicity upper bound is finite.

Section~\ref{S:Upper} discusses upper regularity of
filter coverings (and dually, lower regularity of ideal
coverings). We show in this section that anomalous
filter and ideal coverings cannot be regular. We also
give an example of an atomic filter covering in an
algebraic lattice that is not upper regular, and thus
not regular.

Section~\ref{S:RegFilCov} gives a proof that certain
filter and ideal coverings are regular. In particular, we
show that in a meet-continuous lattice, if the
multiplicity upper bound of a covering $x\prec y$ if
finite, then the corresponding filter covering is
regular.

Section~\ref{S:Chief} discusses the application of
these ideas to generalizing the Jordan-Holder Theorem.

Section~\ref{S:Proof} applies the theory to defining the
steps in the proof of a proposition from given premises.

We will talk almost entirely about
modular lattices, complete in most cases, except in
Section~\ref{S:Proof}, where we will talk about the
distributive lattice underlying a boolean algebra
$\mathbf B$, and complete distributive lattices
constructed from it.

\section*{0.  Preliminaries}

The reader should know about \emph{modular lattices} and
\emph{distributive lattices}, and that distributive
lattices are modular. The reader should also know about
\emph{complete lattices}.

We denote the least element of any lattice, if
one exists, by $\bot$, and the greatest element by $\top$.
If $x\leq y$ are elements of a lattice $L$, then we denote
by $\I_L[x,y]$, or simply $\I[x,y]$, the \emph{interval
sublattice} of elements $z$ such that $x\leq z\leq y$.

A \emph{covering} is a pair $\pair xy$ of elements such
that $x<y$ and $\I[x,y]$ has only $x$ and $y$ as elements.
We say that $x$ is \emph{covered} by $y$, or $x\prec y$.
We will often say $x\prec y$ not only to
state that $x$ is covered by $y$, but also to denote a
pair $\pair xy$ satisfying the covering relation.

If $L$ is a lattice, we say that an element $m\in L$ is
\emph{meet-irreducible} if $x>m$, $y>m$ imply $x\wedge
y>m$. If $L$ is complete, then we say that $m$ is
\emph{strictly meet-irreducible} if for all
$S\subseteq L$ such that $s\in S$ implies $s>m$,
$\bigwedge S>m$. Note that if $m$ is strictly
meet-irreducible, then there is a unique element $m'$
such that $m\prec m'$.

 If
$x$,
$y$,
$z$, and
$w\in L$ with
$x\leq y$ and
$z\leq w$, we write $\pair xy\nearrow\pair zw$ when
$\pair xy$
\emph{transposes up to} $\pair zw$, i.e., when $y\wedge
z=x$ and $y\vee z=w$. When pairs $\pair xy$, $\pair
zw$, such that $x\leq y$ and $z\leq w$, are related by
the symmetric and transitive closure of $\nearrow$, we
say that they are \emph{projectively equivalent}, or
$\pair xy\sim\pair zw$.

Projective equivalence classes of coverings in
modular lattices will be of fundamental importance to
us.  The projective equivalence class of a covering
$x\prec y$ will be denoted by $[x\prec y]$.

% For example, we might say
%$x\prec y\nearrow z\prec w$, if $\pair xy\nearrow\pair
%zw$ and $y$ and $w$ cover $x$ and $z$, respectively.

A lattice $L$ is a \emph{chain} if the
natural ordering in $L$ is a total order. Also, if $L$ is a
lattice, and $C\subseteq L$, then $C$ is called a
\emph{chain in $L$} if in the ordering inherited from $L$,
$C$ is a chain. If $C$ is a chain in $L$, then we say
$C$ is \emph{maximal} if no larger subset of $L$ is a
chain in $L$.

A complete lattice $L$ is \emph{meet-continuous} if for
all
$a\in L$ and $D\subseteq L$ such that $D$ is
\emph{directed upward} (i.e., $d$, $d'\in D$ imply there
exists $d''\in D$ such that $d\leq d''$ and $d'\leq d''$)
we have
\[a\wedge\bigvee D=\bigvee_{d\in D}(a\wedge d).\] A
lattice with the dual property is called
\emph{join-continuous}.

 For some other important concepts of
lattice theory that we shall mention--in particular,
lattices which are \emph{algebraic} or
\emph{coalgebraic}--we refer to texts on lattice theory
such at
\cite{B} and \cite{C-D}. We will use the fact that
algebraic lattices are meet-continuous, and
coalgebraic lattices are join-continuous.

In section~\ref{S:Chief}, we also assume an acquaintance
with the basic concepts of Universal Algebra, as defined,
for example, in \cite{B-S}. In
particular, the concept of a \emph{congruence} will be
used, and that of the \emph{congruence lattice} of an
algebra. The reader should know that the congruence
lattice of an algebra is always algebraic, and hence,
meet-continuous.

The reader should know about cardinal and ordinal
numbers, as used in transfinite induction. If $\kappa$
is a cardinal number, then $\Succ\kappa$ will stand for
the smallest cardinal number strictly greater than
$\kappa$.

\section{Multiplicity and Multiplicity
Bounds}\label{S:Mult}

\subsection{$C$-Multiplicity}
If $L$ is a modular lattice, $C$ is a
chain in $L$, and $u$, $v\in C$ are such that $u\prec
v$, then we say that the covering $u\prec v$ is
\emph{in} $C$. If $x\prec y$ is a covering in $L$, $C$
is a chain in $L$, and the set of coverings
$u\prec v$ in $C$ such that $u\prec
v\sim x\prec y$ has cardinality $n$, then we say that the
\emph{$C$-multiplicity} of $x\prec y$ (in $L$), denoted by
$\mu_C[x\prec y]$, is
$n$.

\subsection{Weak regularity}

We say that a covering $x\prec y$ is \emph{weakly
regular} if $\mu_C[x\prec y]$, for
maximal chains $C$, is a number $\mu[x\prec y]$, the
\emph{multiplicity} of $x\prec y$, independent of $C$.

If $x\prec y$ is weakly regular, with finite
multiplicity, then we can talk not only about the
multiplicity of $x\prec y$ in $L$, but in any interval
sublattice $\I[a,b]$ of $L$ where $a<b$:

\begin{theorem} If $L$ is a modular lattice and $x\prec
y$ is weakly regular in $L$, with finite multiplicity,
then given $a$, $b\in L$ with $a<b$, the number of
coverings $u\prec v$ equivalent to $x\prec y$ in any
maximal chain of elements in the interval sublattice
$\I[a,b]$ is a number $\mu^{a,b}[x\prec y]$ independent
of the particular chain $C$.
\end{theorem}

\begin{proof} Any two maximal chains $C$, $C'\in\I[a,b]$
can be completed to maximal chains in $L$ by including
the elements of the same maximal chains in $\I[\bot,a]$
and $\I[a,\top]$ (where we first adjoin a $\bot$ and
a $\top$ to $L$ if not already present). Then we use the
fact that
$\mu_{\tilde C}[x\prec y]=\mu_{\tilde C'}[x\prec y]$.
\end{proof}

\subsection{Multiplicity upper bounds and lower bounds}

Let $x\prec y$ be a covering in a modular lattice $L$.
We define $\upsilon[x\prec y]$, the \emph{multiplicity
upper bound} of
$x\prec y$ in
$L$, to be the least cardinal number $\nu$ such that for
every chain $C$ in $L$, $\mu_C[x\prec y]<\nu$.
We define $\lambda[x\prec y]$, the \emph{multiplicity
lower bound} of $x\prec y$ in $L$, to be the least
cardinal $\nu$ such that $\mu_C[x\prec y]=\nu$ for some
maximal chain $C$.

\begin{proposition} If $x\prec y$ is weakly regular in
$L$, then
\[\lambda[x\prec y]=\mu[x\prec y]\]
and
\[\upsilon[x\prec y]=\Succ\mu[x\prec y].\]
\end{proposition}

\subsection{Distributive lattices}

The $C$-multiplicity is severely constrained for
distributive lattices:

\begin{theorem}\label{T:DistMult} If $L$ is a distributive
lattice,
$C$ is a maximal chain in $L$, and $x\prec y$ is a
covering in $L$, then $\mu_C[x\prec
y]$ is $0$ or
$1$.
\end{theorem}

\begin{proof}
Assume that $u\prec v\nearrow z\prec w$ and $u'\prec
v'\nearrow z\prec w$. Then we claim that $u\wedge u'\prec
v\wedge v'$ and $u\wedge u'\prec v\wedge v'\nearrow
u\prec v$. For,
\begin{align*}(v\wedge v')\wedge u
&=(v\wedge v')\wedge(v\wedge z)\cr
&=(v\wedge z)\wedge(v'\wedge z)\cr
&=u\wedge u',
\end{align*}
and
\begin{align*}(v\wedge v')\vee u
&=(v\vee u)\wedge(v'\vee u)\cr
&=v\wedge(v'\vee u)\cr
&=v\wedge(v'\vee(v\wedge z))\cr
&=v\wedge(v'\vee v)\wedge(v'\vee z)\cr
&=v\wedge(v'\vee v)\wedge w\cr
&=v.
\end{align*}
Similarly, $u\wedge u'\prec v\wedge v'\nearrow u'\prec
v'$.

It follows that if $u\prec v\sim u'\prec v'$, then
we must have some covering $z\prec w$ such that $z\prec
w\nearrow u\prec v$ and
$z\prec w\nearrow u'\prec v'$. Therefore,
$u$, $v$, $u'$, and $v'$ cannot all be elements of
the same chain $C$.
\end{proof}

\begin{remark} As a result of his theorem, if $L$ is a
distributive lattice, and $a<b\in L$, then we can talk
about the set of weakly regular coverings $x\prec y$ in
$\I[a,b]$. We will do so in the last section of this
paper.
\end{remark}

\section{The Theory of Regular Coverings}\label{S:RegCov}

\subsection{Upper regular and lower regular
coverings}
We say that a covering $x\prec y$ in a complete
modular lattice $L$ is \emph{upper regular} if, for every
chain $I$, and mapping taking elements $i\in I$ to
coverings $x_i\prec y_i$ of $L$, projectively equivalent
to $x\prec y$ and such that $i<j$ implies $x_i\prec
y_i\nearrow x_j\prec y_j$, we have
$\bigvee_ix_i\prec\bigvee_iy_i$ (rather than the only
other possibility, for a modular lattice, which would be
$\bigvee_ix_i=\bigvee_iy_i$.) The property dual to upper
regularity, we call \emph{lower regularity}. We say that a
covering is \emph{regular} if it is both upper regular and
lower regular.

Clearly, whether or not a covering is upper regular,
lower regular, or regular depends only on the projective
equivalence class of the covering.

The importance of the concept of regularity comes from a
generalization of the Dedekind-Birkhoff Theorem,
proved in \cite{R97}:

\begin{theorem}\label{T:MultInv} If $L$ is a complete
modular lattice, and $C$, $C'$ are any two maximal chains
in~$L$, then for every regular covering $x\prec y$,
$\mu_C[x\prec y]=\mu_{C'}[x\prec y]$.
\end{theorem}

Thus, if $x\prec y$ is regular, we can drop the $C$
from $C$-multiplicity and speak of the \emph{multiplicity}
$\mu[x\prec y]$ of
$x\prec y$ in $L$. In other words, if $x\prec y$ is
regular, then $x\prec y$ is weakly regular.

A partial converse to Theorem~\ref{T:MultInv}:

\begin{theorem}\label{T:FiniteRegular}
Let $L$ be a
complete modular lattice. If $x\prec y$ is weakly regular,
and furthermore, $\mu[x\prec y]$ is finite, then $x\prec
y$ is regular.
\end{theorem}

\begin{proof} Let $I$ be a chain, and let coverings
$x_i\prec y_i\sim x\prec y$ be indexed by $I$,
such that
$i<j$ implies $x_i\prec y_i\nearrow x_j\prec y_j$. Then,
for any arbitrary $i\in I$, consider the chain
$x_i\leq\ldots\leq
x_j\leq\ldots\leq\bigvee_ix_i\leq\bigvee_iy_i$ and the
chain $x_i\prec y_i\leq\ldots\leq
y_j\leq\ldots\leq\bigvee_iy_i$. If we take any refinement
of the first of these chains to a maximal chain $C$, we can
find a refinement of the second chain to a maximal
chain $C'$, by letting $C'$ consist of the elements of $C$
less than or equal to $x_i$, the lattice elements $c\vee
y_i$ for $c\in C$ such that $x_i\leq c\leq\bigvee_ix_i$,
and the elements of $C$ greater than or equal to
$\bigvee_iy_i$. By the modular law, the coverings in $C$
between $x_i$ and $\bigvee_ix_i$ correspond in a
one-to-one fashion with the coverings in $C'$ between
$y_i$ and $\bigvee_iy_i$, and corresponding coverings are
projectively equivalent. Since $\mu_C[x\prec
y]=\mu_{C'}[x\prec y]$, and that number is finite,
we must have $\bigvee_ix_i\prec\bigvee_iy_i$, proving that
$x\prec y$ is upper regular. Lower regularity is proved
similarly.
\end{proof}

In order to apply Theorem~\ref{T:MultInv}, it
helps to know which coverings are regular. Some
preliminary observations in this direction are as
follows:  If $L$ is finite, or of finite height, then all
coverings in
$L$ are regular. It is easy to see
that in any complete, modular,
meet-continuous lattice, every covering is upper
regular.  Dually, in any complete,
modular, join-continuous lattice, every covering is lower
regular.

\section{Coverings in Filter and Ideal Lattices}
\label{S:Coverings}

In this section, we will explore coverings in filter and
ideal lattices. 
If $L$ is a lattice, a \emph{filter} in $L$ is a nonempty
subset $F$ such that if $x\in F$, and $y\ge x$, then
$y\in F$, and also, if $x$, $y\in F$, then $x\wedge y\in
F$. If $F$ and $G$ are filters in $L$, we say that
$F\leq G$ if $G\subseteq F$. With this partial
ordering, the filters of a lattice $L$ with $\top$ form a
lattice $\Fil L$ which is complete and coalgebraic. We
have $F\vee G=F\cap G$, while $F\wedge G=\{\,z\in L:z\geq
x\wedge y\hbox{, for some }x\in F\hbox{ and }y\in G\,\}$.
If $S\subseteq L$ is a nonempty subset, we write $\Fg(S)$
for the smallest (in the sense of set inclusion) filter
containing $S$, called \emph{the filter generated by
$S$}. An important special case is $\Fg\{\,x\,\}$, the
\emph{principal filter} generated by $x\in L$, which is
$\{\,y\in L:y\geq x\,\}$. The mapping $x\mapsto
\Fg\{\,x\,\}$ is a lattice homomorphism embedding $L$ into
$\Fil L$. As another important example of a filter, if
$m$ is a meet-irreducible element, then we denote by
$F_{>m}$ the set of elements of $L$ strictly greater than
$m$. $F_{>m}$ is obviously a filter, and is principal iff
$m$ is not just meet-irreducible, but strictly
meet-irreducible.

The dual concept, that of an \emph{ideal}, leads to the
lattice of ideals $\Idl L$, which is complete and
algebraic. If $S\subseteq L$, we write $\Ig S$ for the
smallest ideal containing $S$, and call it \emph{the
ideal generated by $S$}. If $x\in L$, the \emph{principal
ideal} generated by $x$, $\Ig\{\,x\,\}=\{\,y\in L:y\leq
x\,\}$ is an important example. $\Idl L$ is ordered
by inclusion as opposed to $\Fil L$, which is ordered
by reverse inclusion. The mapping
$x\mapsto \Ig\{\,x\,\}$ is a lattice homomorphism
embedding $L$ into~$\Idl L$.

Both the lattices $\Fil L$ and $\Idl L$ satisfy every
lattice-theoretic identity satisfied by $L$; in
particular, they are modular if $L$ is modular.

For the most part, we will concentrate our attention on
filters of a modular lattice $L$, leaving to the reader the
dualization of the statements and proofs of the theorems
to yield similar results about ideals of $L$.

\subsection{Filter coverings $F\prec G$ and $\mathcal
M(F-G)$}

If $L$ is a lattice and $S\subseteq L$, then we denote
by~$\mathcal M(S)$ the set of maximal elements of $S$
(in the partial ordering of $L$). If $F\prec G$ is a
covering in
$\Fil L$, or, as we say, a \emph{filter covering}, we will
be particularly interested in
$\mathcal M(F-G)$. We have

\begin{lemma}\label{T:MaxLemma}
Let $L$ be a lattice, and $F$, $G$, $F'$,
$G'\in\Fil L$ such that $F\prec G$, $F'\prec G'$, and
$F\prec G\nearrow F'\prec G'$. Then
\begin{enumerate}
\item[(1)] $F'-G'=F'\cap(F-G)$, and
\item[(2)] $\mathcal M(F'-G')=F'\cap\mathcal M(F-G)$.
\end{enumerate}
\end{lemma}

\begin{proof} (1): $F'\cap(F-G)=(F'\cap F)-(F'\cap
G)=F'-G'$.

(2): $x\in\mathcal M(F'-G')\implies x\in F'\cap(F-G)$ by
(1). If, in addition, $y>x$, then $y\in G'$ which implies
$y\in G$. Thus, $x\in F'\cap\mathcal M(F-G)$. On the
other hand, if $x\in F'\cap\mathcal M(F-G)$, then $x\in
F'-G'$ by (1), and $y>x\implies y\in G\implies y\in
F'\cap G=G'$. Thus, $x\in\mathcal M(F'-G')$.
\end{proof}

\subsection{Filter coverings and maximal based filters}
If $F\prec G$ is a filter covering in $\Fil L$, then any
$x\in F-G$ determines a principal filter $\Fg\{\,x\,\}$.
We have $F\prec G\nearrow
\Fg\{\,x\,\}\prec(\Fg\{\,x\,\}\vee G)$. Let
$H=\Fg\{\,x\,\}\vee G=\Fg\{\,x\,\}\cap G$. We say that a
pair $\pair xH$, such that $\Fg\{\,x\,\}\prec H$, is a
\emph{maximal based filter} with $x$ as its \emph{base}.
We say that $\pair xH$ is a maximal based filter
\emph{determined by} $F\prec G$.

We now define separate concepts of $\nearrow$ and
projective equivalence for maximal based filters. Suppose
$\pair xG$ and $\pair yH$ are maximal based filters. We say
$\pair xG\nearrow \pair yH$ if $x\le y$, $y\not\in G$, and
$H=\Fg\{\,y\,\}\vee G=\Fg\{\,y\,\}\cap G$. We call the
symmetric, transitive closure of this relation \emph{
projective equivalence of maximal based filters} and again
use the symbol~$\sim$. It is easy to see that
the relation of projective equivalence of maximal based
filters is a subset of the relation of projective
equivalence on filters, restricted to maximal based
filters viewed as filter coverings. That is, $\pair
xG\sim\pair yH$ implies $\Fg\{\,x\,\}\prec G\sim
\Fg\{\,y\,\}\prec H$. The converse is also true:

\begin{lemma}\label{T:ProjEquiv}
Let $L$ be a lattice. Given two filter
coverings $F\prec G$ and $H\prec K$ in $\Fil L$, and
given $x\in F-G$ and $y\in H-K$, $F\prec G\sim H-K$
iff $\pair x{\Fg\{\,x\,\}\vee G}\sim\pair
y{\Fg\{\,y\,\}\vee K}$.
\end{lemma}

\begin{proof} It suffices to prove that if $F\prec
G\nearrow H\prec K$, then for any such $x\in F-G$ and
$y\in H-K$, we have
$\pair x{\Fg\{\,x\,\}\vee G}\sim\pair y{\Fg\{\,y\,\}\vee
K}$.

$y\in H-K$ implies $y\in F-G$, so we have $F=G\wedge
\Fg\{\,y\,\}$.
Thus, there exists $g\in G$ such that $g\wedge y\leq
x$. Then we have $\pair{g\wedge y}{\Fg\{\,g\wedge
y\,\}\vee G}\nearrow\pair x{\Fg\{\,x\,\}\vee G}$ and
$\pair{g\wedge y}{\Fg\{\,g\wedge y\,\}\vee G}\nearrow\pair
y{\Fg\{\,y\,\}\vee K}$, whence
$\pair x{\Fg\{\,x\,\}\vee G}\sim\pair y{\Fg\{\,y\,\}\vee
K}$.
\end{proof}

\begin{corollary}\label{T:ProjEqFil} If $x\prec y$ and
$z\prec w$, then
$x\prec y\sim z\prec w$ in $L$ iff $\Fg\{\,x\,\}\prec
\Fg\{\,y\,\}\sim \Fg\{\,z\,\}\prec \Fg\{\,w\,\}$ in $\Fil
L$.
\end{corollary}

\begin{proof} This follows from Lemma~\ref{T:ProjEquiv}
and from the fact that $x\prec y\nearrow z\prec w$ iff
$\pair x{\Fg\{\,y\,\}}\nearrow\pair z{\Fg\{\,w\,\}}$.
\end{proof}

\subsection{Atomic filter coverings}
Let $\pair xF$ be a maximal based filter.
We say that $\pair xF$, or a filter covering $F'\prec G'$
such that $\Fg\{\,x\,\}\prec F\sim F'\prec G'$, is \emph{
atomic} if $F$ is principal.

As an example, if $m\in L$
is strictly meet-irreducible, then $\pair
m{F_{>m}}$ is an atomic maximal based filter.

\begin{theorem} \label{T:Atomic} Let $L$ be a modular
lattice. The set of atomic maximal based filters in $L$,
and the set of atomic filter coverings, are closed under
projective equivalence. If $F\prec G$ is a filter
covering, and $m\in F-G$ is strictly
meet-irreducible, then $F\prec G$ is atomic.
\end{theorem}

\begin{proof} Let $\pair xF\nearrow\pair yG$. If $F$ is
principal, say $F=\Fg\{\,x'\,\}$, then
$G=\Fg\{\,y\,\}\cap\Fg\{\,x'\,\}=\Fg\{\,y\vee x'\,\}$ is
also principal.

On the other hand, if $G$ is principal, say
$G=\Fg\{\,y'\,\}$, then let $\bar x\in F$ be such that
$x=\bar x\wedge y$, and let $x'=\bar x\wedge y'$. We
cannot have $x'=x$ because both $\bar x$ and $y'$ belong
to $F$. Thus, by modularity, $x'\succ x$. But, this
implies that $F=\Fg\{\,x'\,\}$. Thus, the set of atomic
maximal based filters is closed under projective
equivalence, and by Theorem~\ref{T:ProjEquiv}, the same
is true of the set of atomic filter coverings.

If $F\prec G$ and $m\in F-G$ is strictly
meet-irreducible, then $F\prec G\sim\Fg\{\,m\,\}\prec
(G\cap\Fg\{\,m\,\})$. However, $F_{>m}$ is the unique
cover of $\Fg\{\,m\,\}$ and is principal. It follows that
$G\cap\Fg\{\,m\,\}=F_{>m}$, and
$F\prec G$ is atomic.
\end{proof}

\begin{theorem}\label{T:MAtomic} Let $L$ be a complete,
meet-continuous modular lattice, and let $F\prec G$ be an
atomic filter covering. Then $\mathcal M(F-G)$ is
nonempty and consists of strictly meet-irreducible
elements.
\end{theorem}

\begin{proof} Let $x\in F-G$. Then $F\prec
G\nearrow\Fg\{\,x\,\}\prec(\Fg\{\,x\,\}\vee
G)=\Fg\{\,x'\,\}$ where $x'\succ x$, because $F\prec G$
is atomic. The set of elements $y$ such that $y\geq x$
and
$y\wedge x'=x$ is closed under joins of chains, by
meet-continuity. Then by Zorn's Lemma, $\mathcal
M(\Fg\{\,x\,\}-\Fg\{\,x'\,\})
=\Fg\{\,x\,\}\cap\mathcal M(F-G)$ is nonempty. Thus,
$\mathcal M(F-G)$ is nonempty. It is easy to see that
because $F\prec G$ is atomic, $\mathcal M(F-G)$ consists
of strictly meet-irreducible elements.
\end{proof}

\subsection{Quasi-atomic filter coverings}
We say that a maximal based filter $\pair xF$, or a
filter covering $F'\prec G'$ such that
$\Fg\{\,x\,\}\prec F\sim F'\prec G'$, is
\emph{quasi-atomic} if $F$ is not principal, but
contains an element
$y$ such that
$x<z\leq y$ implies $z\in F$.

As an example, if $m\in L$ is meet-irreducible, but not
strictly meet-irreducible, then $\pair
m{F_{>m}}$ is a quasi-atomic maximal based filter.

\begin{theorem}\label{T:QA} Let $L$ be a modular lattice.
The set of quasi-atomic maximal based filters, and the
set of quasi-atomic filter coverings, are closed under
projective equivalence. If $F\prec G$ is a filter
covering, and
$m\in F-G$ is meet-irreducible but not
strictly meet-irreducible, then
$F\prec G$ is quasi-atomic.
\end{theorem}

\begin{proof} Let $\pair xF\nearrow\pair yG$. If $\pair
xF$ is quasi-atomic, then $F$ contains an element $x'$
such that
$x<z\leq x'\implies z\in F$. Consider the element
$y'=x'\vee y\in G=\Fg\{\,y\,\}\cap F$. If $y<z\leq y'$,
then $x'\wedge z>x$, because by modularity,
$y\vee (x'\wedge z) = (y\vee x')\wedge z=y'\wedge
z=z>y$. Thus, $y<z\leq y'\implies z\in
G=F\cap\Fg\{\,y\,\}$, and
$\pair yG$ is quasi-atomic because if it were atomic,
then $\pair xF$ would also be atomic by
Theorem~\ref{T:Atomic}.

On the other hand, if $\pair yG$ is quasi-atomic, then
there is an element $y'\in G$ such that $y<z\leq y'$
implies $z\in G$. Let $\bar x\in F$ be such that $x=\bar
x\wedge y$, and let $x'=\bar x\wedge y'$. We have $x'\in
F$, so $y\vee x'\in G$. If $x<z\leq x'$, then by
modularity, $z=z\vee(y\wedge x')=(z\vee y)\wedge x'$.
But, $y<z\vee y\leq y'$ because if $y=z\wedge y$, then
$z=(z\vee y)\wedge x'=x$. Thus, $z\vee y\in G$ and $z\in
F$. It follows that $\pair xF$ is atomic or quasi-atomic,
but $\pair xF$ cannot be atomic, because then $\pair yG$
would also be atomic by Theorem~\ref{T:Atomic}.

Now, if $F\prec G$ and $m\in F-G$ is meet-irreducible, but
not strictly meet-irreducible, we have $F\prec
G\nearrow\Fg\{\,m\,\}\prec(\Fg\{\,m\,\}\vee G)$. However,
$F_{>m}$ is the unique cover of $\Fg\{\,m\,\}$ in~$\Fil
L$. Thus, $F\prec G\sim\Fg\{\,m\,\}\prec F_{>m}$, which
is quasi-atomic.
\end{proof}

\begin{theorem}\label{T:MQA} Let $L$ be a complete,
meet-continuous modular lattice. If $F\prec G$ is a
filter covering in $L$ that is quasi-atomic, then
$\mathcal M(F-G)$ is a nonempty set of
elements of $L$ that are meet-irreducible, but not
strictly meet-irreducible.
\end{theorem}

\begin{proof} Similar to the proof of
Theorem~\ref{T:MAtomic}. Instead of $\Fg\{\,x\,\}\vee
G=\Fg\{\,x'\,\}$, we have an $x'\in\Fg\{\,x\,\}\vee G$
such that $\Fg\{\,x\,\}\vee G=\Fg\{\,z\mid x<z\leq
x'\,\}$. The set of elements $y$ such that $y\geq x$ and
$y\wedge x'=x$ is again closed under joins of chains by
meet-continuity, and nonempty by Zorn's Lemma. Thus
$\mathcal M(\Fg\{\,x\,\}-(\Fg\{\,x\,\}\vee G))$ is
nonempty, and so is $\mathcal M(F-G)$ by
Lemma~\ref{T:MaxLemma}.
\end{proof}

\subsection{Anomalous filter coverings}

We say that a maximal based filter $\pair xF$, or a
filter covering~$F'\prec G'$ such that $\Fg\{\,x\,\}\prec
F\sim F'\prec G'$, is
\emph{anomalous} if it is neither atomic nor
quasi-atomic.

Recall that $x\in L$ is called \emph{finitely
decomposable} if $x$ is a finite meet of
meet-irreducible elements. For an example of an anomalous
filter covering, let
$x\in L$ be an element which is not finitely
decomposable. (This is possible only if
$L$ does not satisfy  the ascending chain condition.) Let
$G$ be the filter generated by the set of
finitely decomposable elements of $L$ that
are greater than $x$. (This is the same as
the filter generated by the set of meet-irreducible
elements of
$L$ that are greater than
$x$.) We have $x\notin G$ because otherwise, $x$ would
be finitely decomposable. By Zorn's Lemma, there is a
filter
$F\leq G$ such that $\Fg\{\,x\,\}<F\leq G$ and $F$ is
minimal (in the ordering of $\Fil L$) for that property.
Then by Theorems~\ref{T:MAtomic} and \ref{T:MQA}, $\pair
xF$ is an anomalous maximal based filter, because it
cannot be atomic or quasi-atomic. The following theorem
shows, among other things, that this example is typical:

\begin{theorem}\label{T:Anom} The set of anomalous
maximal based filters, and the set of anomalous filter
coverings, are closed under projective equivalence. If
$F\prec G$ is an anomalous filter covering, then
$\mathcal M(F\prec G)$
 is empty, and $F-G$ contains no elements that are
finitely decomposable.
\end{theorem}

\begin{proof} The sets of aomic and quasi-atomic filter
coverings are closed under projective equivalence. Since
the set of anomalous filter coverings comprises the rest
of the filter coverings, it is also closed under
projective equivalence. If $F\prec G$ and $x\in\mathcal
M(F-G)$, then clearly $x$ is meet-irreducible. Thus, if
$F\prec G$ is anomalous, $\mathcal M(F-G)$ must be empty
by Theorems \ref{T:MAtomic}~and~\ref{T:MQA}. Finally, if
$x$ is finitely decomposable, then
$x=\bigwedge_{i=1}^nm_i$ where the $m_i$ are
meet-irreducible. If $x\in F-G$, then $m_i$ also belongs
to $F-G$ for some $i$, because $G$ is closed under
finite meets. This would imply that $F\prec G$ was atomic
or quasi-atomic. 
\end{proof}

\subsection{A counterexample}
In working with the $\nearrow$ relation and filters, we
might make the following conjecture:
 Let $L$ be a modular
lattice, and $F$, $G$, $H$, $K$ filters such that
$F\prec G$, $H\prec K$, and $F\prec G\nearrow H\prec
K$. If $x\in F-G$, then there exists $w\in H-K$ such
that $x\leq w$.
However, this is false: 

\begin{example}
Consider the modular lattice
known as $M_5$, with its elements labeled as follows:
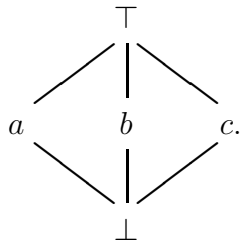
\[
\thicklines
\begin{picture}(100,100)
\put(10,50){$a$}
\put(52,50){$b$}
\put(90,50){$c.$}
\put(50,90){$\top$}
\put(50,10){$\bot$}
\put(55,84){\line(0,-1){20}}
\put(55,24){\line(0,1){20}}
\put(50,24){\line(-4,3){30}}
\put(59,24){\line(4,3){30}}
\put(50,84){\line(-4,-3){30}}
\put(59,84){\line(4,-3){30}}
\end{picture}
\]
Let
$F=M_5$, $G=\Fg\{\,b\,\}=\{\,b,\top\,\}$,
$H=\Fg\{\,c\,\}=\{\,c,\top\,\}$,
and $K=\Fg\{\,\top\,\}=\{\,\top\,\}$. Then $F\prec
G\nearrow H\prec K$. Observe that we have
$a\in F-G$, but no element $w\in H-K$ such that $a\leq w$.
\end{example}

\section{Stability Theorems}
\label{S:Stability}

We will shortly begin to address the
question of when a covering in $\Fil L$ is regular. 
First, however, we pose and answer some other important
questions, such as, under what circumstances do regular
coverings in $L$ remain regular after the embedding from
$L$ into
$\Fil L$ or
$\Idl L$? We also examine stability of multiplicity, in
case a covering is regular, and of the multiplicity
upper bound.  We continue to focus on
$\Fil L$. A lemma:

\begin{lemma} \label{T:LowTech} Let $L$ be a complete
modular lattice, and
$x\prec y$ a covering in $L$ which is lower regular. If
$F\prec G$ is a covering in $\Fil L$
and $F\prec G\sim \Fg\{\,x\,\}\prec \Fg\{\,y\,\}$,
then let $f=\bigwedge F$ and $g=\bigwedge G$; we have
\begin{enumerate}
\item[(1)]
$f\prec g$,
\item[(2)]
$f\prec g\sim x\prec y$, and
\item[(3)]
$\mathcal
M(\Fg\{\,f\,\}-\Fg\{\,g\,\})=\mathcal
M(F-G)$.
\end{enumerate}
\end{lemma}

\begin{proof} For all $w\in F-G$, $\Fg\{\,w\,\}\wedge
G=F$. It follows that $w\wedge g=f$ and hence,
$\Fg\{\,w\,\}\wedge \Fg\{\,g\,\}=\Fg\{\,f\,\}$.

We must show that $f\neq g$, which we
will show by showing that if $w\in F-G$, then we cannot
have $g\leq w$, or in other words, we
cannot have $\bigwedge_\nu g_\nu\leq w$ for any
$\kappa$-tuple
$\{\,g_\nu\,\}_{\nu<\kappa}$ of elements of $G$, for any
cardinal number $\kappa$, where $\nu$ runs through
ordinals less than $\kappa$. Assume the contrary, where
$\kappa$ is the least cardinal possible. By modularity,
and the fact that $F\prec G$ is atomic, we can assume
w.l.o.g.\ that
$g_\nu\wedge w\prec g_\nu$ for each $\nu$. For each
ordinal~$\nu\leq \kappa$, let
$h_\nu=\bigwedge_{\nu'<\nu}g_{\nu'}$. By the minimality of
$\kappa$, we have $h_\nu\not\leq w$ if $\nu<\kappa$, so we
must have
$h_\nu\wedge w\prec h_\nu$ for $\nu<\kappa$. Note $\kappa$
cannot be finite, because
$G$ is a filter. By lower regularity, we have
$h_\kappa\wedge w\prec h_\kappa$, contradicting the
assumption that
$h_\kappa\leq w$.

We have $\Fg\{\,g\,\}\not\leq F$, $\Fg\{\,g\,\}\leq G$, and
$\Fg\{\,f\,\}=\Fg\{\,g\,\}\wedge F$, whence
$\Fg\{\,f\,\}\prec \Fg\{\,g\,\}$, proving
(1). Also, $\Fg\{\,f\,\}\prec
\Fg\{\,g\,\}\nearrow \Fg\{\,w\,\}\prec \Fg\{\,w\,\}\cap
G$ for any $w\in F-G$. (2) follows by
corollary~\ref{T:ProjEqFil}.

We have $\mathcal M(F-G)\subseteq\mathcal
M(\Fg\{\,f\,\}-\Fg\{\,g\,\})$. For, let $m\in\mathcal
M(F-G)$. Then $m\in\Fg\{\,f\,\}$, and we cannot have
$m\in\Fg\{\,g\,\}$, because then we would have $g\leq m$,
contradicting the fact that $m\wedge g=f$.
$m\in\mathcal M(\Fg\{\,f\,\}-\Fg\{\,g\,\})$ because $m$
is strictly meet-irreducible.

 On the other hand, let
$m\in\mathcal M(\Fg\{\,f\,\}-\Fg\{\,g\,\})$. Since
$m\notin\Fg\{\,g\,\}$, $m\notin G$. It suffices to
show $m\in F$, because, $m$ being strictly
meet-irreducible, $M\in F-G$ will imply $m\in\mathcal
M(F-G)$. Let
$\kappa$ be the least cardinal number such that some
$\kappa$-tuple
$\{\,f_\nu\,\}_{\nu<\kappa}$ of elements of $F$, where
$\nu$ runs through ordinals less than $\kappa$, satisfies
$\bigwedge_\nu f_\nu\leq m$. $m$ is strictly
meet-irreducible because $\Fg\{\,f\,\}\prec \Fg\{\,g\,\}$
is atomic. Let $m'$ be the unique cover of $m$,
and for each
$\nu\leq\kappa$, define
$u_\nu=\bigwedge_{\nu'<\nu}f_{\nu'}$. We have
$m\vee u_\nu\geq m'$ if $\nu<\kappa$, by the minimality
of~$\kappa$. Thus, for each $\nu<\kappa$, we have by
modularity
\[m\wedge
u_\nu\prec m'\wedge u_\nu
\nearrow m\prec m'.\]
If $\nu'<\nu<\kappa$, then it is easy to see that
\[m\wedge u_\nu\prec m'\wedge u_\nu\nearrow
m\wedge u_{\nu'}\prec m'\wedge u_{\nu'}.\]
Now, if $\kappa$ is infinite, then by the lower regularity
of
$x\prec y$, we must have $m\wedge u_\kappa\prec m'\wedge
u_\kappa$. However, this is absurd because $u_\kappa\leq
m$. Thus,
$\kappa$ is finite. It follows that $m\in F$, proving
(3).
\end{proof}

\begin{theorem}\label{T:RegStable}
 Let $L$ be a complete, meet-continuous modular
lattice, and
$x\prec y$ a covering in $L$ which is regular. Then
$\Fg\{\,x\,\}\prec \Fg\{\,y\,\}$ is regular in $\Fil L$.
\end{theorem}

\begin{proof} It suffices to prove upper regularity,
because $\Fil L$ is coalgebraic, which implies that all
coverings are automatically lower regular.
Suppose given
$F_i\prec G_i\sim
\Fg\{\,x\,\}\prec
\Fg\{\,y\,\}$, indexed by $i\in I$, for some chain $I$, and
such that $F_i\prec G_i\nearrow F_j\prec G_j$ for $i<
j$. We must show that $\bigvee_iF_i<\bigvee_iG_i$.

For each $i\in I$, let $M_i=\mathcal M(F_i-G_i)$. We
have $M_j\subseteq M_i$ for $i<j$ by
Lemma~\ref{T:MaxLemma}, and
$\bigcap_iM_i\subseteq\bigvee_iF_i-\bigvee_iG_i$. For,
if $x\in M_i$ for all $i$, then $x\in F_i$ for all $i$ so
$x\in\bigvee_iF_i=\bigcap_iF_i$, and $x\notin G_i$ for
all $i$, so $x\notin\bigvee_iG_i=\bigcap_iG_i$. Thus, it
suffices to show that
$\bigcap_iM_i$ is nonempty.

For each $i$, let $f_i=\bigwedge F_i$ and
$g_i=\bigwedge G_i$. Since $x\prec y$ is lower regular, we
have
$f_i\prec g_i$ and
$f_i\prec g_i\sim x\prec y$ by Lemma~\ref{T:LowTech}(1)
and (2), and for
$i< j$, we have $f_i\prec g_i\nearrow f_j\prec g_j$.
For, $f_i\leq f_j$, $f_i\prec g_i$, and
$g_i\not\leq f_j$, because for any $x\in F_j-G_j$, we
have $F_i=G_i\wedge\Fg\{\,x\,\}$ and consequently,
$f_i=g_i\wedge x$. It follows that
$f_i=g_i\wedge f_j$. We also have $g_i\leq g_j$ and
$f_j\prec g_j$, whence $g_j=g_i\vee f_j$.

Since $x\prec y$ is upper regular, we have
$f\prec g$ where $f=\bigvee_if_i$ and $g=\bigvee_ig_i$.
This implies that $\mathcal M(\Fg\{\,f\,\}-\Fg\{\,g\,\})$
is nonempty, by Theorem~\ref{T:MAtomic}.
Let $x\in\mathcal M(\Fg\{\,f\,\}-\Fg\{\,g\,\})$. Then
$x\geq f_i$ for all $i$ so $x\in\Fg\{\,f_i\,\}$ for all
$i$.
On the other hand, $x\not\geq g$ so there is an $i$ such
that $x\notin\Fg\{\,g_i\,\}$. If $j> i$ then $g_j\geq
g_i$, so
$x\notin\Fg\{\,g_j\,\}$. If $y>x$, then $y\geq g$ and
$y\in G_k$ for all $k$. Thus, $x\in\mathcal
M(\Fg\{\,f_j\,\}-\Fg\{\,g_j\,\})$ for all $j\geq i$.
On the other hand, if
$j<i$, then by Lemma~\ref{T:MaxLemma},
$\mathcal
M(\Fg\{\,f_i\,\}-\Fg\{\,g_i\,\})=\Fg\{\,f_i\,\}\cap\mathcal
M(\Fg\{\,f_j\,\}-\Fg\{\,g_j\,\})$, implying that
$x\in\mathcal M(\Fg\{\,f_j\,\}-\Fg\{\,g_j\,\})$ in this
case as well. It follows that
$\bigcap_i\mathcal M(\Fg\{\,f_i\,\}-\Fg\{\,g_i\,\})$
is nonempty. However, by Lemma~\ref{T:LowTech}(3),
$\mathcal M(\Fg\{\,f_i\,\}-\Fg\{\,g_i\,\})=M_i$ for each
$i$. Thus,
$\bigcap_iM_i$ is nonempty, and $F\prec G$ is regular.
\end{proof}

Now, we consider the multiplicity:

\begin{lemma}\label{T:MultStableLemma} If $L$ is a
complete modular lattice and
$x\prec y$ is lower regular in $L$, then for any maximal
chain $C\subseteq L$ and any maximal chain
$\bar C\subseteq\Fil L$ refining the image of $C$ in $\Fil
L$, $\mu_{\bar
C}[\Fg\{\,x\,\}\prec\Fg\{\,y\,\}]=\mu_C[x\prec y]$.
\end{lemma}

\begin{proof} Let $C$ be a maximal chain in $L$, and
$\bar C$ a maximal chain in $\Fil L$ refining the image
of~$C$.

If we have a principal filter $\Fg\{\,u\,\}\in\bar C$,
then we must have $u\in C$. For, if $u\notin\bar C$,
then there must exist $c\in C$ such that $u$ and $c$ are
not comparable. But, $\Fg\{\,c\,\}\in\bar C$, so either
$\Fg\{\,c\,\}<\Fg\{\,u\,\}$, implying $c<u$, or
$\Fg\{\,u\,\}<\Fg\{\,c\,\}$, implying $u<c$.

If $F$, $G\in\bar C$ are principal and such that
$F\prec G$, say $F=\Fg\{\,u\,\}$ and $G=\Fg\{\,v\,\}$,
then we must have $u$, $v\in C$ with $u\prec v$, and
if $F\prec G\sim\Fg\{\,x\,\}\prec\Fg\{\,y\,\}$ then
$u\prec v\sim x\prec y$.

Let
$F$,
$G\in\bar C$ be such that
$F\prec G$ and it is not true that $F$ and $G$ are both
principal. We will show that $F\prec G$ is not
projectively equivalent to $\Fg\{\,x\,\}\prec
\Fg\{\,y\,\}$.

If $F$ is principal and $G$ is not, then $F\prec G$ is
not of atomic type, is not projectively equivalent to
$\Fg\{\,x\,\}\prec \Fg\{\,y\,\}$, and does not count in the
multiplicity of $x\prec y$.

The case $G$ principal and $F$ non-principal cannot
occur, because $F\prec G$.

The only remaining case is that both $F$ and $G$ are
non-principal. We can assume that $F\prec G$ is atomic
and lower regular, since otherwise $F\prec
G\sim\Fg\{\,x\,\}\prec\Fg\{\,y\,\}$ is impossible. If we
had
$c\in C$ with
$c\in F-G$, then we would have $F<\Fg\{\,c\,\}<G$; thus,
we must have
$C\cap F=C\cap G$ in order to have $F\prec G$. Denote
this set by
$D$. Since $F\prec G$ is atomic, let $z\in F-G$
and $w\in G$ with $z\prec w$. For each $d\in D$, we have
$d\wedge z\prec d\wedge w$, and for $d<d'\in D$, we have
$d\wedge z\prec d\wedge w\nearrow d'\wedge z\prec d'\wedge
w$. Since $z\prec w$ is lower regular in $L$, we would
have
\[(\bigwedge D)\wedge z=\bigwedge_{d\in D} d\wedge
z\prec\bigwedge_{d\in D} d\wedge w=(\bigwedge D)\wedge
w.\]
Now, if we had $\bigwedge D\in C-D$, we would have
$\Fg\{\,\bigwedge D\,\}<F\prec G$, implying that
$\bigwedge D\wedge z=\bigwedge D\wedge w$, which is
impossible.

The only other possibility is $\bigwedge D\in D$. In
this case, we claim that we must have $\bigwedge D\leq
w$ and $\bigwedge D\wedge z\in C$. For, if $c\in  C-D$,
then $\Fg\{\,c\,\}\leq F$, implying that $c\leq\bigwedge
D\wedge z$. Then because $C$ is maximal, we must have
$\bigwedge D\wedge z\in C-D$ and $\bigwedge D\wedge
w=\bigwedge D\in D$. However, this contradicts the fact
that $F\cap C=D$, because $\bigwedge D\wedge z\in F$. It
follows that the case $F\prec G$ atomic, lower regular,
and neither $F$ nor $G$ principal is impossible.
\end{proof}

It follows from the Lemma that we have

\begin{theorem}\label{T:MultStable} If $x\prec y$ is
regular in
$L$, then $\mu[\Fg\{\,x\,\}\prec\Fg\{\,y\,\}]=\mu[x\prec
y]$.
\end{theorem}

Finally, we examine the stability properties
of the multiplicity upper bound.

\begin{theorem}\label{T:MUB} If $x\prec y$ is a covering
in a modular lattice $L$, then $\upsilon[x\prec
y]\leq\upsilon[\Fg\{\,x\,\}\prec\Fg\{\,y\,\}]$, with
equality if
$\upsilon[x\prec y]$ is finite or countable. In any case,
$\upsilon[x\prec y]$ infinite implies
$\upsilon[\Fg\{\,x\,\}\prec\Fg\{\,g\,\}]$ infinite.
\end{theorem}

\begin{proof}
Let $C$ be a chain in $L$. Then  the image of $C$ in
$\Fil L$ can be refined to a maximal chain $\bar C$ in
$\Fil L$, and it is clear that $\mu_C[x\prec
y]\leq\mu_{\bar C}[\Fg\{\,x\,\}\prec\Fg\{\,y\,\}]$.

Now, let
\[F_1\prec G_1\leq F_2\prec G_2\leq\ldots\leq
F_n\prec G_n
\]
where $F_i\prec G_i\sim\Fg\{\,x\,\}\prec \Fg\{\,y\,\}$ for
all $i$. Then $\exists u_i\in F_i-G_i$, $v_i\in G_i$ such
that $u_i\prec v_i$ and $u_i\prec v_i\sim x\prec y$. For
each $i$, define
$c_i=\bigwedge_{j\geq i}u_i$,
$d_i=v_i\wedge\bigwedge_{j>i}u_i$.
Then $c_i=d_i\wedge u_i$, $d_i\in G_i$, $c_i\in F_i-G_i$,
and $d_i\leq v_i$, implying that $c_i\prec d_i\sim
x\prec y$ and
\[c_1\prec d_1\leq c_2\prec d_2\leq\ldots
\leq c_n\prec d_n.
\]
It follows that $\upsilon[x\prec y]\geq n$, and combined
with the fact that $\upsilon[x\prec
y]\leq\upsilon[\Fg\{\,x\,\}\prec\Fg\{\,y\,\}]$, this
implies that
$\upsilon[x\prec
y]=\upsilon[\Fg\{\,x\,\}\prec\Fg\{\,y\,\}]$ if
$\upsilon[x\prec y]$ is finite or countable, and that
$\upsilon[x\prec y]$ is infinite if
$\upsilon[\Fg\{\,x\,\}\prec\Fg\{\,y\,\}]$ is.
\end{proof}

\section{$\upsilon[x\prec y]$,
$\lambda[\Fg\{\,x\,\}\prec\Fg\{\,y\,\}]$, and
$\lambda[\Ig\{\,x\,\}\prec\Ig\{\,y\,\}]$}
\label{S:Relationship}

In this section, we consider the relationship
between the multiplicity upper bound of a covering, and
the multiplicity lower bound of the corresponding filter
covering or ideal covering. As usual, we focus on
filter coverings, leaving the dual result to be stated
by the reader.

Consider the function $\Lambda:\mathbb N\to\mathbb N$,
where $\mathbb N$ stands for the natural numbers, defined
recursively as follows:
\[\Lambda(n) = \begin{cases}
0, &n=0\\
1+\Lambda(\lfloor\sqrt n\rfloor-1),&n>0.
\end{cases}
\]

\begin{lemma}\label{T:Lambda} We have
\begin{enumerate}
\item[(1)] $\Lambda(n)\geq 0$ for all $n$
\item[(2)] $\Lambda$ is increasing; i.e.,
$n<n'\implies\Lambda(n)\leq\Lambda(n')$
\item[(3)] $\lim_{n\to\infty}\Lambda(n)=\infty$
\item[(4)] $\Lambda(n)\leq\lfloor\sqrt n\rfloor$ for all
$n$.
\end{enumerate}
\end{lemma}

\begin{proof} (1) is clear.

To prove (2), note that we have $\Lambda(0)=0$ and
$\Lambda(1)=1$, so (2) is true for $n<n'\leq 1$. If (2)
is true for $n<n'\leq\bar n>1$ then for $n\leq\bar n+1$,
$\lfloor\sqrt n\rfloor - 1\leq\bar n$, and the
square root function is also increasing, whence
$\Lambda(n')-\Lambda(n)=\Lambda(\lfloor\sqrt{n'}\rfloor-1)-
\Lambda(\lfloor\sqrt n\rfloor-1)\geq 0$ if $n<n'\leq\bar
n+1$. Thus, (2) follows by induction.

To prove (3), we use (2) and note that if $\Lambda(n)=m$,
then $\Lambda((n+1)^2)=m+1$.

A computation shows that the inequality (4) holds for all
$n\leq 10$. Suppose (4) holds for $n\leq\bar n>10$, and
let us prove it is true for $n=\bar n+1$. We have
by the induction hypothesis
\[\Lambda(\bar n+1)=1+\Lambda(\lfloor\sqrt{\bar
n+1}\rfloor-1)\leq 1+\sqrt{\sqrt{\bar n+1}-1}.\]
Squaring, we have
\[\Lambda(\bar n+1)^2=1+2\sqrt{\sqrt{\bar
n+1}-1}+\sqrt{\bar n+1}-1\leq 3\sqrt{\bar n+1}\leq\bar
n+1.\]
Thus, (4) holds for $n=\bar n+1$, and by induction, for
all $n$.
\end{proof}

\begin{theorem} Let $L$ be a complete, meet-continuous
modular lattice, and $x\prec y$ a covering in $L$. Let
\[b_1\succ a_1\geq b_2\succ a_2\geq\ldots
\geq b_n\succ a_n,\]
where $x\prec y\sim a_i\prec b_i$ for all $i$. If $\bar
C$ is a maximal chain in $\Fil L$, then $\mu_{\bar
C}[\Fg\{\,x\,\}\prec\Fg\{\,y\,\}]\geq \Lambda(n)$.
\end{theorem}

\begin{proof}
Let $m_1\in L$ be maximal for the property that $m_1\geq
a_1$ but $m_1\not\geq b_1$. (By meet-continuity, the set
of such elements is closed under joins of chains, so
a maximal such element exists by Zorn's Lemma.)
Then $m_1$ is strictly meet-irreducible, and has a unique
cover, $m_1'=m_1\vee b_1$.

Now, let $m_2\in L$ be maximal for the property that
$m_2\geq a_2$, $m_2\not\geq b_2$, and $m_2\leq m_1$. We
cannot have $m_2=m_1$, because $m_1\geq a_1\geq b_2$.
$m'_2=m_2\vee b_2$ is the unique cover of $m_2$ in the
interval $\I[\bot,m_1]$, because if $x>m_2$ and $x\leq
m_1$, we must have $x\geq b_2$.

Similarly, we successively choose $m_3$, $\ldots$, $m_n$
such that $m_i$ is maximal among elements $x$ such that
$x\geq a_i$, $x\not\geq b_i$, and $x\leq m_{i-1}$, and we
obtain covers $m'_i=m_i\vee b_i$. We have
\[m'_1\succ m_1\geq m'_2\succ m_2\geq\ldots\geq m'_n\succ
m_n.\]

For each $i$, let $F_i$ be the join of all elements of
$\bar C$ containing $m_i$, and $G_i$ the meet of all
elements of $\bar C$ not containing $m_i$. We have
$F_i$, $G_i\in\bar C$ because the maximal chain $\bar
C$ is be closed under joins and meets. Clearly,
$F_i\prec G_i$ for all $i$.

The mapping $i\mapsto F_i\prec G_i$ sends each $i$ to the
unique covering in $\bar C$ such that $m_i\in F_i-G_i$,
and thus partitions the ordered set $\{\,1,\ldots,n\,\}$
into intervals. If $\{\,i,\,i+1,\,\ldots,j\,\}$ is one of
these intervals, then we claim that $m'_i\in G_i$. For,
if $i=1$ then $m_i$ is strictly meet-irreducible, and so
we must have $m'_i\in G_1$. If $i>1$, then we have
$m_{i-1}\in F_{i-1}$, so $m_{i-1}\in G_i$, because,
$\bar C$ being a chain,
$G_i\leq F_{i-1}$. If $m'_i$ did
not belong to
$G_i$, then there would be an element~$g\in G_i$ such
that
$m'_i\wedge g=m_i$. Then we would have
$m'_i\wedge(g\wedge m_{i-1})=m_i$, but $m_i\neq g\wedge
m_{i-1}$ because $G_i$ is closed under meets and does not
contain $m_i$. This is impossible, because $m_i$ is
meet-irreducible in $\I[\bot,m_{i-1}]$. Thus, the claim
that $m'_i\in G_i$ is proved. It follows that $F_i\prec
G_i\sim\Fg\{\,x\,\}\prec\Fg\{\,y\,\}$.

Clearly, we have $m'_{i+1}$, $m_{i+1}$, $\ldots$, $m'_j$,
$m_j\in F_i-G_i$, as well as $m_i$. Thus, since we have
shown that $F_i-G_i$ is atomic, $m_{i+1}$, $\ldots$,
$m_j$ have unique covers $\bar m_{i+1}$, $\ldots$, $\bar
m_j\in G_i$. Defining $\bar m'_k=\bar m_k\vee m'_k$ for
$k=i+1$, $\ldots$, $j$, we obtain
\[\bar m'_{i+1}\succ\bar m_{i+1}\geq\ldots\geq\bar
m'_j\succ\bar m_j\]
and each $\bar m_k\prec\bar m'_k\sim x\prec y$.

Now, either the number of intervals is $\geq\lfloor\sqrt
n\rfloor$, or the cardinality of the largest interval is
$\geq\lfloor\sqrt n\rfloor$. In the first case, we have
$\mu_{\bar
C}[\Fg\{\,x\,\}\prec\Fg\{\,y\,\}]\geq\Lambda(n)$ by
Lemma~\ref{T:Lambda}(4). In the second case, by induction
on $n$ (and noting that recursive application of the
construction of this proof will find filter coverings
above $F_i\prec G_i$ and therefore distinct from it), we
also have
$\mu_{\bar C}[\Fg\{\,x\,\}\prec\Fg\{\,y\,\}]\geq
1+\Lambda(\sqrt n-1)=\Lambda(n)$.
\end{proof}

\begin{corollary}\label{T:InfStable} Let $L$ be a
complete, meet-continuous modular lattice. If
$\upsilon[x\prec y]$ is infinite, then so is
$\lambda[\Fg\{\,x\,\}\prec\Fg\{\,y\,\}]$.
\end{corollary}

\section{Upper Regularity of Filter Coverings and Joins of
Chains}\label{S:Upper}

In this section, we consider the issue of upper regularity
of filter coverings, and show anomalous filter coverings
cannot be regular, because they cannot be upper regular.
We also give an example of an atomic filter covering, in
an algebraic lattice, which is not upper regular, showing
that upper regularity alone of
$x\prec y$ does not imply upper regularity of
$\Fg\{\,x\,\}\prec\Fg\{\,y\,\}$, even if the lattice is
meet-continuous.

Filter coverings are always lower regular, because $\Fil
L$ is coalgebraic, thus join-continuous. Thus, if a
filter covering is upper regular, it must be regular.

A necessary condition for a filter
covering to be upper regular is easy to state:

\begin{theorem}\label{T:JoinsChains} Let $L$ be a complete
modular lattice, and $F$, $G\in\Fil L$ such that $F\prec
G$. If
$F\prec G$ is upper regular, then $F\prec G$ is
closed under joins of chains.
\end{theorem}

\begin{proof} Suppose $C$ is a chain in $F-G$. For each
$c\in C$, define $F_c=\Fg\{\,c\,\}$ and $G_c=G\vee
F_c=G\cap F_c$. Then for all $c\in C$, $F\prec G\nearrow
F_c\prec G_c$, and if
$c$,
$c'\in C$ with
$c\leq c'$, we have $F_c\prec G_c\nearrow F_{c'}\prec
G_{c'}$. Since
$F\prec G$ is upper regular,
$\bigvee_cF_c\prec\bigvee_cG_c$. However,
$\bigvee_cF_c=\Fg\{\,\bigvee C\,\}$. If $\bigvee C\in
G$, then we would have $\bigvee_cF_c=\bigvee_cG_c$. Thus,
$\bigvee C\in F-G$.
\end{proof}

\begin{corollary}\label{T:AnomNotReg} If $L$ is a complete
modular lattice, $F$, $G\in\Fil L$ with $F\prec G$, and
$F\prec G$ is anomalous, then $F\prec G$ is not
upper regular.
\end{corollary}

\begin{proof} By the Theorem, if $F\prec G$ is upper
regular, then
$F-G$ is closed under joins of chains. Then, by Zorn's
Lemma,
$F-G$ has maximal elements. However, this is impossible
for anomalous $F\prec G$ by Theorem~\ref{T:Anom}.
\end{proof}

Some sufficient conditions for the preceding necessary
condition to hold:

\begin{theorem}\label{T:JChainX} Let $L$ be a complete,
modular, meet-continuous lattice, and $F\prec G$ a
covering in
$\Fil L$ which is atomic or quasi-atomic.
Then $F-G$ is closed under joins of chains.
\end{theorem}

\begin{proof} Use meet-continuity
as in the proof of Theorem~\ref{T:MAtomic} or
Theorem~\ref{T:MQA}.
\end{proof}

\begin{theorem}\label{T:JChainY} Let $L$ be a complete
modular lattice, and $x\prec y$ a covering in $L$ which
is upper regular. If $H\prec K\sim \Fg\{\,x\,\}\prec
\Fg\{\,y\,\}$, then $H-K$ is closed under joins of chains.
\end{theorem}

\begin{proof} Let $C$ be a chain in $H-K$, and let $c\in
C$. We have $\Fg\{\,c\,\}\vee K=\Fg\{\,q\,\}$ for some
$q\in K$ such that $c\prec q$ and $c\prec q\sim x\prec y$.
For each $c'$, $c''\in C$ such that $c\leq c'\leq c''$, we
have $c'\prec q\vee c'\nearrow c''\prec q\vee c''$. Then
$\bigvee C=\bigvee_{c'\geq c}c'\prec\bigvee_{c'\geq
c}(c'\vee q)$ by the upper regularity of $x\prec y$. If
we had $\bigvee C\in K$, then we would have
$c=q\wedge\bigvee C\in K$, which is absurd. It follows that
$\bigvee C\in H-K$.
\end{proof}

There follows an example of a meet-continuous lattice,
having an atomic filter covering which is
not regular:

\begin{example}\label{E:Vec} Let $V$ be the
infinite-dimensional real vector space of sequences of
real numbers, only a finite number of which are nonzero.
Let
$L$ be the lattice of subspaces of $V$. For
each finite set $S\subseteq\mathbb N$ of cardinality
$\geq 2$, consider the subspace
\[A_S=\{\,\langle a_0, a_1,\ldots\rangle\in V :
s\in S\implies a_s=0\,\},\]
and the subspace
\[B_S=\{\,\langle b_0,b_1,\ldots\rangle\in V :
s,s'\in S\implies b_s=b_{s'}\,\}.\]
Note that $A_S\prec B_S$ for each $S$, and if
$S\subseteq S'$ then $A_{S'}\prec B_{S'}\nearrow
A_S\prec B_S$.

For each
$n$, let
\[U_n=\{\,A_S:s\in
S\implies s\geq n\,\},\] and
\[U'_n=\{\,B_S:s\in S\implies s\geq n\,\};\]
the sets $U_n$, $U'_n$ are bases for filters $F_n=\Fg
U_n$ and $G_n=\Fg U'_n$.

We have
$F_n\prec G_n$ for all $n$. For, if $H\in F_n-G_n$, then
$A_S\subseteq H$ for some $S$ such that $s\in S\implies
s\geq n$, but there does not exist an $S'$ such that
$s\in S'\implies s>n$ and $B_{S'}\subseteq H$. In
particular, $B_S\not\subseteq H$.
Then $B_S\subseteq H\vee\mathbb Rv$
where $\mathbb R$ is the field of real numbers, and $v$
has
$1$ in positions
$s$ such that $s\in S$, and
$0$ elsewhere. $H\prec H\vee\mathbb Rv$, because $\mathbb
Rv$ is an atom of~$L$.

We have $A_S\prec B_S\nearrow H\prec(H\vee B_S$. Now let
$A_{\hat S}$ be a basic element of $F_n$, and we will
show that $H\cap B_{\bar S}\leq A_{\hat S}$ for some
$\bar S$. If we had $A_{\hat S}\in G_n$ already, this
would be trivial. If $A_{\hat S}\in F_n-G_n$, however,
we have $B _{\hat S}\in G_n$. Then let $\bar S=S\cup\hat
S$. We have $A_{\bar S}=B_{\bar S}\cap H=B_{\bar S}\cap
B_S\cap H=B_{\bar S}\cap A_S=A_{\bar S}\leq A_{\hat S}$.
Thus, $F_n\prec G_n$. This argument has also shown that
the covering $F_n\prec G_n$ is atomic.

Also, if $n'>n$ then we claim that $F_n\prec
G_n\nearrow F_{n'}\prec G_{n'}$. To prove this, it
suffices to prove
$F_n\prec G_n\nearrow F_{n+1}\prec G_{n+1}$. We have
$F_n<F_{n+1}$ and
$G_n<G_{n+1}$, and if $S$ is such that $s\in S\implies
s>n$, and $A_S\in F_n$, then either $A_S\in F_{n+1}$, or
$n+1\in S$. We may
assume that $\card S>2$ since the $A_S$ for such $S$ and
such that $A_S\in F_n$ form a base for $F_n$. Then
$A_S=B_{\{\,n+1,j\,\}}\cap A_{S-\{\,n+1\,\}}$, where
$j\neq n+1$ is any other element of $S$. Thus,
$F_n=G_n\wedge F_{n+1}$. On the other hand, if $H\in
G_{n+1}$, then there is an $S$ such that $s\in S\implies
s>n+1$ and $B_S\subseteq H$. $B_S\in G_n$ and $B_S\in
F_{n+1}$, since $A_S\subseteq B_S$, so
$B_S\in G_n\vee F_{n+1}$. Thus, $G_n\vee F_{n+1}\leq
G_{n+1}$. The claim follows.

We already proved that
$F_n\prec G_n$ is atomic for each $n$. However, $\bigvee
F_n=\bigvee G_n=\{\,V\,\}$. Thus, the coverings $F_n\prec
G_n$ are not upper regular, and so are not regular.
\end{example}

\section{Regularity of
Filter Coverings} \label{S:RegFilCov}

\begin{lemma} Let $L$ be a complete lattice and
$S\subseteq L$, where $S\neq\emptyset$. Then the
following are equivalent:
\begin{enumerate}
\item[(1)] $\Fg S$ is principal, and
\item[(2)] $\bigwedge S\in\Fg S$.
\end{enumerate}
\end{lemma}

Let $L$ be a complete, meet-continuous modular lattice. If
$F\prec G$ is an atomic or quasi-atomic filter covering,
such that the equivalent conditions of the Lemma are
satisfied, then we say that
$F\prec G$ is
\emph{principally bounded}.

\begin{theorem} \label{T:FDReg} Let $L$ be a complete,
meet-continuous modular lattice, and $F\prec G$
a covering in $\Fil L$ which is atomic or
quasi-atomic. If every filter covering $H\prec K$
such that
$H\prec K\sim F\prec G$ is principally bounded, then
$F\prec G$ is regular.
\end{theorem}

\begin{proof} Let $I$ be a chain, and $F_i\prec G_i$ be
filter coverings projectively equivalent to $F\prec G$
and such that $i<j$ implies $F_i\prec G_i\nearrow
F_j\prec G_j$. For each $i$, let $q_i=\bigwedge\mathcal
M(F_i-G_i)$, and let $q=\bigvee_iq_i$. We have $q_i\in
F_i$ for each $i$, so $q\in\bigvee_iF_i$.

On the other hand, we have $q\notin\bigvee_iG_i$. For,
$q_i\leq$ some element of $\mathcal M(F_i-G_i)$, whence
$q_i\notin G_i$. Also, $G_j\geq G_i$ for $j>i$, so
$q_i\notin G_j$ in that case.

If $j<i$, then $G_i=G_j\vee F_i=G_j\cap F_i$. However,
$q_i\in F_i$. Thus,
$q_i\notin G_j$.

So, $q_i\notin G_j$ for all $i$ and $j$.
Now,
$F_j-G_j$ is closed under joins of chains, by
Theorem~\ref{T:JChainX}. Thus,
$q=\bigvee_iq_i\notin G_j$ for all $j$. It follows that
$q\notin\bigvee_iG_i$.

Thus, $q\in\bigvee_iF_i-\bigvee_iG_i$, implying that
$\bigvee_iF_i\prec \bigvee_iG_i$, and that $F\prec G$ is
upper regular, hence regular.
\end{proof}

\begin{corollary}\label{T:DistReg} Let $L$ be a complete,
meet-continuous distributive lattice. If
$F\prec G$ is a covering in $\Fil L$ which is atomic
or quasi-atomic, then $F\prec G$ is regular.
\end{corollary}

\begin{proof}
If $H\prec K\sim F\prec G$, then $H\prec K$ is not
anomalous, so
$\mathcal M(H-K)$ is nonempty by Theorem~\ref{T:MAtomic}
and Theorem~\ref{T:MQA}. Furthermore, by distributivity,
it has cardinality one, proving that
$H\prec K$ is principally bounded. The Corollary then
follows from Theorem~\ref{T:FDReg}.
\end{proof}

\begin{theorem} \label{T:UpsFinite} Let $L$ be a complete,
meet-continuous modular lattice. If
$x\prec y$ is a covering in
$L$ such that
$\upsilon[x\prec y]$ is finite, then
$[\Fg\{\,x\,\},\Fg\{\,y\,\}]$ is regular, with
multiplicity equal to~$\upsilon[x\prec y]$.
\end{theorem}

\begin{proof} If $\Fg\{\,x\,\}\prec\Fg\{\,y\,\}$ were not
principally bounded, $\mathcal
M(\Fg\{\,x\,\}-\Fg\{\,y\,\})$ would contain a sequence of
strictly meet-irreducible elements $m_1$, $m_2$, $\ldots$
such that for all $n$,
$\bigwedge_{i<n}m_i\neq\bigwedge_{i\leq n}m_i$. For
each $i$, let
$m'_i$ be the unique cover of $m_i$. Then if
$y\in L$ is such that $y\not\leq m_i$,
we have $y\wedge m_i\prec y\wedge m'_i$. For, we must 
have $y\vee m_i\geq m'_i$, whence $y\vee m'_i=y\vee m_i$.
If we had $y\wedge m_i=y\wedge m'_i$, then the elements
$y$, $m_i$, $m'_i$, $y\vee m_i$, and $y\wedge m_i$ would
form a sublattice isomorphic to the lattice $N_5$, which
cannot happen in a modular lattice.

It follows that
\[m_1\wedge m'_2\succ m_1\wedge m_2\geq
m_1\wedge m_2\wedge m'_3\succ m_1\wedge m_2\wedge
m_3\geq\ldots,\]
with $\bigwedge_{i\leq
n}m_i\prec(\bigwedge_{i<n}m_i)\wedge m'_n$, the
general covering in the chain, equivalent to $x\prec
y$ for all
$i$. This sequence of elements can be refined to a
maximal chain
$C$ such that $\mu_C[x\prec y]$ is infinite. However,
this is contrary to the assumption that $\upsilon[x\prec
y]$ is finite.

Thus, $\Fg\{\,x\,\}\prec\Fg\{\,y\,\}$ is principally
bounded, and regular by Theorem~\ref{T:FDReg}.
\end{proof}

\section{Lattice-theoretic Chief Factors}\label{S:Chief}

Suppose we have an algebra $A$ (in the sense of
universal algebra) which has a modular congruence
lattice. Then we can apply the preceding theory to the
congruence lattice $\Con A$ and talk about the chief
factors of $A$, obtaining a generalization of the nice
multiplicity result seen in the Jordan-Holder Theorem.

\subsection{Coverings of rank $\mathcal F$ and
lattice-theoretic chief factors of rank $\mathcal F$}
If we have a regular covering $x\prec y$ in the
lattice
$\mathcal F(L)$, where the functor
$\mathcal F$ is some composite of the functors $\Fil$ and
$\Idl$, then we say that $x\prec y$ is a \emph{covering of
$L$ of rank $\mathcal F$}. Then, if $L=\Con A$, we
say that a covering in $L$ of rank $\mathcal F$ is a
\emph{ lattice-theoretic chief factor} of $A$ \emph{of
rank
$\mathcal F$}.

In the theories of finite groups and finite-height
modules, where the lattices involved have finite height,
it is standard practice to assign a group or module to a
covering, obtaining a \emph{chief
factor}, or, in case
$A$ is a module, a \emph{composition factor}, of $A$. In
order to do something similar for an arbitrary
congruence-modular algebra $A$, it is necessary to assign
some type of algebraic object to each lattice-theoretic
chief factor. We leave to future investigations the
question of the manner in which this may be done
generally. (We have taken some small steps toward such a
theory in \cite{R92},
\cite{R96}, and \cite{R98}.) However, the lattice-theoretic
chief factors themselves are of interest,
because their multiplicities are invariants of the algebra.

Thus, in the remainder of this section, unless otherwise
specified, $L$ will denote the lattice~$\Con A$, for
some algebra $A$ such that $\Con A$ is modular. Since
$\Con A$ is algebraic, we also are assuming that $L$ is
meet-continuous.

\subsection{The case when $L=\Con A$ satisfies the
descending chain condition}
If $L$ satisfies the
descending chain condition, then all coverings in $L$
are lower regular, all filters are principal, and all
filter coverings are of atomic type. In fact, $\Fil
L\cong L$. Since $L$ is meet-continuous,
coverings in $L$ (and $\Fil L$) are also upper regular.
Thus, in this situation, all coverings in $L\cong \Fil
L$ are regular.
This result was stated but not proved in \cite{R97};
it must be admitted, however, that as an
example for the application of the ideas in that
paper, and of lattice-theoretic chief
factors, it is vacuous.

\subsection{The case when $L=\Con A$ is distributive}
A better example presents itself when $L$ is
distributive. Then, Corollary~\ref{T:DistReg} shows that
every filter covering $F\prec G$ of atomic or strictly
quasi-atomic type is regular. We do not know how many
coverings there may be in $\Fil L$ which are not of
anomalous type.
 However, we note that $\Fil L$ is provided
with a profuse supply of coverings, by which we
mean, somewhat informally, that if $\alpha$, $\beta\in L$
with $\alpha<\beta$, then there is at least one
covering in $\Fil L$ between $\Fg\{\,\alpha\,\}$ and
$\Fg\{\,\beta\,\}$. (This is easy to prove using Zorn's
Lemma.) We can then apply the dual of
Corollary~\ref{T:DistReg} to
$\Idl\Fil L$, because $\Fil L$ is coalgebraic. The
conclusion is that there is a profuse supply of regular
coverings in~$\Idl\Fil L$, i.e., a profuse supply of
lattice-theoretic chief factors of $A$ of rank
$\Idl\circ\Fil$. We have regularized
coverings in $\Fil L$ by the embedding into $\Idl\Fil L$.

Those coverings in $\Idl\Fil L$ which are not known
to be regular can be regularized by considering
them in $\Fil\Idl\Fil L$, where they become regular
by Corollary~\ref{T:DistReg}, as long as $L$ is
distributive. And so on. By this method, any covering in
$\mathcal F(L)$, for any functor $\mathcal F$ which is a
nonempty composite of
$\Idl$ and $\Fil$, can be regularized by applying either
$\Fil$, or $\Idl$, and similarly, any covering in any
distributive lattice whatever can be regularized by
applying either $\Idl\circ\Fil$, or $\Fil\circ\Idl$.

Because of Theorem~\ref{T:DistMult}, all of the regular
coverings that arise in this way have multiplicity one.

\subsection{Modular but not distributive $L=\Con A$}

In
this case, multiplicities higher than $1$ are possible.
If $x\prec y$ is a covering in $L$, such that
$\upsilon[x\prec y]$ is finite, then
$\Fg\{\,x\,\}\prec\Fg\{\,y\,\}$ is regular by
Theorem~\ref{T:UpsFinite}. Thus, the embedding from $L$
into $\Fil L$ regularizes such coverings, and the
multiplicity of the corresponding filter covering is
$\upsilon[x\prec y]$ by Theorem~\ref{T:MUB}.

On the other
hand, if
$\upsilon[x\prec y]$ is infinite, then by
Corollary~\ref{T:InfStable}, so is
$\lambda[\Fg\{\,x\,\}\prec\Fg\{\,y\,\}]$. Thus, any
maximal chain $\bar C$ in $\Fil L$ has an infinite number
of coverings in it that are equivalent~to
$\Fg\{\,x\,\}\prec\Fg\{\,y\,\}$, and we can say that the
multiplicity of the filter covering is infinite, even
though we may not be able to say that that multiplicity
is a well-defined cardinal number.

\section{The Steps in the
Proof of a Theorem}\label{S:Proof}

\newcommand\True{{\scriptstyle\mathbf T}}
\newcommand\False{{\scriptstyle\mathbf F}}

Another application of these ideas is a method of
formalizing the steps necessary and sufficient to prove a
given proposition from given premises. Our
treatment of this will use a simple
Logic framework.

Suppose we have a set  $\mathbf P$ of ``propositions,''
which can in principal be determined to either be true,
or not. Then we have two truth values $\True$ and
$\False$, and for any proposition $P$ we can say that
$\mathcal T(P)$ (the truth value of $P$) takes
values $\True$ and
$\False$.  Given any $n$ propositions $P_1$, $\ldots$,
$P_n$, and any $n$-ary function $f$ with arguments
consisting of truth values, we can formulate a new,
synthetic proposition $f(\vec P)$ with truth value
$f(\mathcal T(P_1),\ldots,\mathcal T(P_n))$. If we
consider the truth values as elements of the two-element
boolean algebra $\{\,\True,\False\,\}$, then the
functions obtainable by compositions of the ordinary
logical connectives give us this, because the
two-element boolean algebra has the property of being
\emph{primal}--i.e., the property that every finitary
function can be constructed from the basic operations.
We will use the symbols
$\wedge$, $\vee$,
$\lnot$, $\rightarrow$ with their usual meanings, along
with $\True$ and $\False$. In fact, it is convenient to
replace our original set of propositions $\mathbf P$ by a
boolean algebra $\mathbf B$ free on $\mathbf P$ as set of
generators. (Or, if $\mathbf P$ already has some or all of
the logical connectives, by a quotient of such a free
boolean algebra.) The assignment $\mathcal T$ of truth
values can then be extended to a boolean algebra
homomorphism from
$\mathbf B$ to the two-element
boolean algebra. Henceforth,
\emph{proposition} shall mean an element of $\mathbf B$.

We will write $\top$ and $\bot$ for the maximum and
minimum elements of $\mathbf B$. The \emph{underlying
lattice} of $\mathbf B$ is just $\mathbf B$, forgetting
the unary operation $\lnot$. A \emph{filter} of $\mathbf
B$ is the same as a filter of the underlying lattice.

 Given some sort of calculus of proving
propositions, consisting of finitary rules of
inference, we assume that the rules of inference include
a small set of \emph{trivial} rules of inference, and
otherwise we call them \emph{nontrivial} rules of
inference. The trivial rules of inference are the rule
that we can infer $P\wedge Q$ from $P$ and $Q$, for any
elements $P$,
$Q\in\mathbf B$, and the rule that for any $P$,
$Q\in\mathbf B$, if $P\leq Q$, then we can infer $Q$
from $P$. Note that \emph{modus ponens}, the rule that
we can infer $Q$ from $P$ and $P\to Q$ (or $\lnot P\vee
Q$) will thus be considered a trivial rule of inference,
because $P\wedge(\lnot P\vee Q)=(P\wedge\lnot
P)\vee(P\wedge Q)=P\wedge Q\leq Q$.

Note that if the ordering of $\mathbf B$ provides that
$P\leq Q$ whenever $Q$ can be proved from $P$, then the
second trivial rule of inference would actually subsume
all the rules of inference, rendering our analysis of the
situation vacuous. Thus, we want to consider a situation
where $\mathbf B$ does not have such an ordering.

 We say
$P\vdash Q$ ($S\vdash Q$, where $S\subseteq\mathbf B$) if
$Q$ can be proved from $P$ (from elements of $S$) using
both trivial and nontrivial rules of inference.

\begin{theorem}\label{T:Pretheory} If $T\subseteq\mathbf
B$, then the following are equivalent:
\begin{enumerate}
\item
$T$ is a filter of the underlying lattice of $\mathbf B$,
and
\item
$T$ is closed under application of the trivial rules of
inference, and contains $\top$.
\end{enumerate}
\end{theorem}

We call a set $T$, satisfying the equivalent conditions
of Theorem~\ref{T:Pretheory}, a \emph{pretheory}. A
\emph{theory} is quite often defined as a set of
propositions closed under the rules of inference. We use
the term \emph{pretheory} to suggest that a filter in
$\mathbf B$ is a forerunner of a theory, and will not
have occasion to mention theories further.

Now, $\mathbf B$ may or may not be complete or algebraic
when viewed as a lattice, and may or may not have any
coverings at all, but $\Fil\mathbf B$ is
complete and coalgebraic, and has a profuse supply of
coverings. Thus, by the dual of corollary~\ref{T:DistReg},
$\Idl\Fil\mathbf B$ has a profuse supply of regular
coverings. (We use the word \emph{profuse} in the
informal sense of the previous section.) We call regular
coverings in $\mathcal F(\Fil\mathbf B)$, where
$\mathcal F$ is a functor as used in
section~\ref{S:Chief}, \emph{steps of order
$\mathcal F$}. That is, an step of order
$\mathcal F$ is a regular covering in $\mathbf B$ of rank
$\mathcal F\circ\Fil$.  For example, if $T$ and $T'$ are
pretheories such that $T\prec T'$,
$\Ig\{\,T\,\}\prec\Ig\{\,T'\,\}$ is regular (by the dual
of Corollary~\ref{T:DistReg}) and is an step of
order
$\Idl$.

If $T$ is a
pretheory, $P$ is a proposition,
and $\mathcal F$ is a functor as before, then we call
the steps of order $\mathcal F$ in
$\I[\phi(T)\wedge
\phi(\Fg\{\,P\,\}),\phi(T)]$, where
$\phi:\Fil\mathbf B\to\mathcal F\Fil\mathbf B$ is the
natural embedding, the \emph{steps of order
$\mathcal F$ in the proof of $P$ from $T$.} Note that we
do not assume
$T\vdash P$.

If we
have an instance of a rule of inference which infers $Q$
from $P_1$, $\ldots,\,P_n$, then we say that that instance
\emph{covers} the set of steps of order
$\mathcal F$ which occur in
$\I[\phi(\Fg\{\,Q\wedge\bigwedge_iP_i\,\}),\phi
(\Fg\{\,\bigwedge_iP_i\,\})]$. If we have a set $\mathbf
N$ of instances of rules of inference, then we say that
$\mathbf N$ \emph{covers} the union of the sets of steps
covered by the individual instances
$N\in\mathbf N$. We also say that $\mathbf N$ covers any
smaller set of steps.

Recall that an \emph{ultrafilter} is a cover of $\bot$ in
$\Fil\mathbf B$.

\begin{theorem}\label{T:EssSteps} Let
$\Ig\{\,T_1\,\}\prec\Ig\{\,T_2\,\}$ be an step of
order $\Idl$, where $T_1\prec T_2$ (i.e., an step
of order $\Idl$, of atomic type), and let $P\in T_1-T_2$.
Then $T_2\wedge\Fg\{\,\lnot P\,\}$ is an ultrafilter of
$\mathbf B$, and $\bot\prec T_2\wedge\Fg\{\,\lnot
P\,\}\nearrow T_1\prec T_2$. This is a one-one
correspondence of ultrafilters with projective equivalence
classes of steps of order $\Idl$, and the steps of order
$\Idl$ in the proof of $P$, from a pretheory $T$,
correspond to the ultrafilters that contain
$T$ but not $P$.
\end{theorem}

\begin{proof} If $T_2\wedge\Fg\{\,\lnot P\,\}$ were
$\bot$, there would be $Q\in T_2$ such that
$Q\wedge\lnot P=\bot$. Then we would have to have $Q\leq
P$, implying that $P\in T_2$ which is not true. Thus,
$T_2\wedge\Fg\{\,\lnot P\,\}$ is an atom, i.e., an
ultrafilter. The ultrafilters $U$ of $\mathbf B$ all
form coverings $\bot\prec U$ which determine distinct
projective equivalence classes, because given two
ultrafilters $U$ and $U'$, if we had $\bot\prec
U\sim\bot\prec U'$, the multiplicity of $\bot\prec U$ in
$\I_{\Fil\mathbf B}[\bot,U\vee U']$ would be $2$, and
this is impossible.
\end{proof}

\begin{theorem}\label{T:EssEmpty} The set of
steps (of any order $\mathcal F$) covered by an instance
of a trivial rule of inference is empty.
\end{theorem}

\begin{theorem}\label{T:Covers} If $T$ is a pretheory and
$P$ is a proposition, then any proof of $P$ from $T$
covers the steps (of any order $\mathcal F$) in
the proof of
$P$ from $T$.
\end{theorem}

\begin{proof} 
Let $N_i$, $i=1$, $\ldots$, $n$ be the
instances of rules of inference in a proof, in order. Let
pretheories $T_i$, $i=0$, $\ldots$, $n$ be defined by
$T_0=T$, $T_i=T_{i-1}\wedge \Fg\{\,Q_i\,\}$ for $0<i\leq
n$, where $Q_i$ is the conclusion of $N_i$. For each $i>0$,
let the set of steps of order $\mathcal F$ in
$\I[\phi(T_i),\phi(T_{i-1})]$ be $E_i$, and the set of
steps of order $\mathcal F$ covered by $N_i$, by
$E'_i$.

We have $E_i\subseteq E'_i$. For, if $N_i$ infers
$Q_i$ from the finite set of propositions $S_i$, we have
\[\pair{T_i}{T_{i-1}}\nearrow
\pair{\Fg\{\,\bar Q_i\,\}}{\Fg\{\,\bar
Q_i\,\}\vee T_{i-1}},\]
where $\bar Q_i=Q_i\wedge\bigwedge S_i$.
However,
$\Fg\{\,\bar Q_i\,\}\vee T_{i-1}\leq
\Fg\{\,\bigwedge S_i\,\}$, because, the $n$-tuple
$\langle N_1,\ldots,N_n\rangle$ being a proof,
$S_i\subseteq T_{i-1}$.

Thus, $\bigcup_iE_i\subseteq\bigcup_iE'_i$,
but the left side is the set of
steps of order $\mathcal F$ in the proof of $P$ from $T$,
and the right side is the set of steps of order
$\mathcal F$ covered by the proof.
\end{proof}

Let $T$, $T'$ be pretheories such that $T'\leq T$, and
let $\mathbf N$ be a set of
instances of rules of inference. For each $N\in\mathbf N$,
let $S_N$ be the (finite) set of premises of $N$, and
$Q_N$ the conclusion. If $T'$ is the join (intersection)
of all pretheories $\bar T\leq T$ such that $N\in\mathbf
N$ and
$S_N\subseteq\bar T$ imply $Q_N\in\bar T$, then we say that
$\mathbf N$ \emph{generates $T'$ from $T$}. In this case,
$T'$ consists of all propositions provable from $T$ using
the elements of $\mathbf N$ as the only instances of
nontrivial rules of inference:

\begin{theorem}\label{T:RefinedToProof} Let $T$, $T'$ be
pretheories with
$T'\leq T$, and let $\mathbf N$ be a set
of instances of rules of inference which generates $T'$
from $T$. Then $T'$ is the set of propositions $P$ such
that there is a finite sequence of elements of $\mathbf N$
that can be refined to a proof of $P$ from $T$ by adding
instances of trivial rules of inference.
\end{theorem}

\begin{proof}
Let $\tilde T$ be that set of propositios, and we will
show that $T'=\tilde T$. Since $T'$ is generated from
$T$ by $\mathbf N$, $T'$ is the intersection (join) of
all pretheories $\bar T\leq T$ such that $N\in\mathbf N$
and $S_N\subseteq\bar T$ imply $Q_N\in\bar T$.

Clearly, $N\in\mathbf N$ and $S_N\subseteq\tilde T$
imply $Q_N\in\tilde T$, because we can construct a proof
of $Q_N$ from proofs of the elemtns of $S_N$. Thus,
$\tilde T\leq T'$.

On the other hand, suppose that $\bar T\leq T$ is such
that $N\in\mathbf N$ and $S_N\subseteq\bar T$ imply
$Q_N\in\bar T$, and let $P\in\tilde T$. The existence of
a proof of $P$ from $T$ using instances from $\mathbf N$
implies that $P\in\bar T$. Thus, $\bar T\leq\tilde T$,
so $T'\leq\tilde T$.

Thus, $T'=\tilde T$.
\end{proof}

Finally, a theorem which shows that covering the steps in
the proof of $P$ from $T$ is not only necessary, but
sufficient:

\begin{theorem} Given pretheories $T$, $T'$
such that $T'\leq T$, a set $\mathbf N$ of
instances of rules of inference that generates $T'$
from $T$, and a proposition $P$, then we have $P\in T'$ iff
$\mathbf N$ covers the steps in the proof of $P$
from $T$.
\end{theorem}

\begin{proof}
If $P\in T'$, then the conclusion follows from
Theorems \ref{T:Covers} and~\ref{T:RefinedToProof}.

If $P\notin T'$, then we have
\[\pair{\Fg\{\,P\,\}\wedge
T'}{T'}\nearrow\pair{\Fg(P)\wedge T}{(\Fg\{\,P\,\}\wedge
T)\vee T'}\] and we have $(\Fg\{\,P\,\}\wedge T)\vee
T'\leq T$. Thus, the steps in the interval
$\I[\Fg\{\,P\,\}\wedge T',T']$ are a subset of the set
of steps in the proof of $P$ from $T$. By Zorn's Lemma,
there is an ideal $J\in\Idl\Fil\mathbf B$ such that
$\Ig\{\,\Fg\{\,P\,\}\wedge T'\,\}<J\prec\Ig\{\,T'\,\}$.
The covering $\Fg\{\,J\,\}\prec\Fg\{\,\Ig\{\,T'\,\}\,\}$
is an step (of order $\Fil\circ\Idl$) in the
proof of $P$ from $T$ that is not covered by $\mathbf N$.
\end{proof}


\begin{thebibliography}{9}

\bibitem{B}
Garrett Birkhoff,
\emph{Lattice Theory (3d Edition)},
AMS Colloquium Publications, Volume XXV,
American Mathematical Society,
Providence,
1967.

\bibitem{B-S}
Stanley Burris and H. P. Sankappanavar,
\emph{A Course in Universal Algebra},
Graduate Texts in Mathematics 78,
Springer--Verlag,
New York,
1981.

\bibitem{C-D}
Peter Crawley and Robert P.\ Dilworth,
\emph{Algebraic Theory of Lattices},
Prentice-Hall,
Inc., Englewood Cliffs, New Jersey,
1973.

%\ref\key MacL\by Saunders Mac Lane
%\book Categories for the Working Mathematician
%\bookinfo Graduate Texts in Mathematics 5
%\publ Springer--Verlag
%\yr 1971
%\endref

\bibitem{R92}
William H.\ Rowan,
\emph{Enveloping Ringoids of Universal Algebras},
Dissertation, University of California
at Berkeley, University Microfilms International,
Ann Arbor, 1992.

\bibitem{R96}
\bysame,
\emph{Enveloping ringoids},
Algebra Universalis
\textbf{35} (1996),
202--229.

\bibitem{R97}
\bysame,
\emph{Regular coverings in complete modular lattices},
Algebra Universalis \textbf{37} (1997),
77--80.


\bibitem{R98}
\bysame,
\emph{The category of directed systems in a category},
Appl. Categorical Structures
\textbf{6} (1998), no.\ 1,
63--86.

\end{thebibliography}
\end{document}